\documentclass[review]{elsarticle}

\usepackage[utf8]{inputenc}
\usepackage{amsmath,amssymb,amsthm,amsfonts}
\usepackage{stackengine}
\usepackage{tcolorbox}
\usepackage{algorithm}
\usepackage{breqn}
\usepackage{stmaryrd}
\usepackage{media9}
\usepackage{caption}
\usepackage{bm}
\usepackage{xcolor}
\usepackage{graphicx}
\graphicspath{{figures/}}
\usepackage{subfig}
\usepackage{epstopdf}
\usepackage{mathtools}
\usepackage{relsize}

\usepackage{algpseudocode}
\algnewcommand{\IfThenElse}[3]{
\State \algorithmicif\ #1\ \algorithmicthen\ #2\ \algorithmicelse\ #3}

\usepackage{natbib}
\usepackage{graphicx}

\usepackage{setspace}

\allowdisplaybreaks

\newcommand{\eugenio}[1]{\textcolor{black}{#1}}
\newcommand{\jon}[1]{\textcolor{black}{#1}}

\newcommand{\Li}{{\,\mathrm{Li}}}
\newcommand{\LLi}{{\,\mathrm{limLi}}}
\newcommand{\U}{{\,\mathrm{U}}}
\newcommand{\LSI}{{\,\mathrm{LSI}\:\!}}
\newcommand{\SQI}{{\,\mathrm{SQI}\:\!}}
\newcommand{\CBI}{{\,\mathrm{CBI}\:\!}}
\newcommand{\HCI}{{\,\mathrm{HCI}\:\!}}

\newcommand{\HCIA}{{\,\mathrm{HCI_A}\:\!}}
\newcommand{\HCIB}{{\,\mathrm{HCI_B}\:\!}}
\newcommand{\HCIC}{{\,\mathrm{HCI_C}\:\!}}

\newcommand{\TRI}{{\,\mathrm{TRI}\:\!}}

\newcommand{\TRIA}{{\,\mathrm{TRI_A}\:\!}}
\newcommand{\TRIB}{{\,\mathrm{TRI_B}\:\!}}
\newcommand{\TRIBR}{{\,\mathrm{TRI_{BR}}\:\!}}
\newcommand{\TRIC}{{\,\mathrm{TRI_C}\:\!}}

\newcommand{\TTI}{{\,\mathrm{TTI}\:\!}}
\newcommand{\TTIA}{{\,\mathrm{TTI_A}\:\!}}
\newcommand{\TTIB}{{\,\mathrm{TTI_B}\:\!}}
\newcommand{\TTIC}{{\,\mathrm{TTI_C}\:\!}}

\newcommand{\PRI}{{\,\mathrm{PRI}\:\!}}
\newcommand{\PRIA}{{\,\mathrm{PRI_A}\:\!}}
\newcommand{\PRIB}{{\,\mathrm{PRI_B}\:\!}}
\newcommand{\PRIC}{{\,\mathrm{PRI_C}\:\!}}

\newcommand{\mints}{\mathbb{N}}
\newcommand{\mC}{\mathbb{C}}
\newcommand{\mZ}{\mathbb{Z}}
\newcommand{\mR}{\mathbb{R}}

\newcommand{\timeseq}{\mathrel{*}=}
\newcommand{\pluseq}{\mathrel{+}=}

\newtheorem{lemma}{Lemma}[section]

\newtheorem{proposition}{Proposition}[section]
\newtheorem{definition}{Definition}[section]

\newtheorem{remark}{Remark}[section]

\begin{document}

\begin{frontmatter}

\title{Exact Subdomain and Embedded Interface Polynomial Integration in Finite Elements with Planar Cuts}


\author{Eugenio Aulisa\corref{mycorrespondingauthor}}
\cortext[mycorrespondingauthor]{Corresponding author}
\ead{eugenio.aulisa@ttu.edu}

\address{Department of Mathematics and Statistics, Texas Tech University, Lubbock TX 79410, USA}
\author{Jonathon Loftin}
\ead{jonathonloftin@saumag.edu}
\address{Department of Mathematics and Computer Science, Southern Arkansas University, Magnolia, AR, 71753, USA}

\begin{abstract}

The implementation of discontinuous functions occurs in many of today's state-of-the-art partial differential equation solvers. However, in finite element methods, this poses an inherent difficulty: \jon{efficient quadrature rules available when integrating functions whose discontinuity falls in the element's interior are for low order degree polynomials, not easily extended to higher order degree polynomials, and cover a restricted set of geometries.} Many approaches to this issue have been developed in recent years. Among them one of the most elegant and versatile is the equivalent polynomial technique. This method replaces the discontinuous function with a polynomial, allowing integration to occur over the entire domain rather than integrating over complex subdomains. Although eliminating the issues involved with discontinuous function integration, the equivalent polynomial tactic introduces its problems. \eugenio{The exact subdomain integration requires a machinery that quickly grows in complexity when increasing the polynomial degree and the geometry dimension, restricting its applicability to lower order degree finite element families.
The current work eliminates this issue. We provide algebraic expressions to exactly evaluate the subdomain integral of any degree polynomial on parent finite element shapes cut by a planar interface. These formulas also apply to the exact evaluation of the embedded interface integral. We provide recursive algorithms that avoid overflow in computer arithmetic for standard finite element geometries: triangle, square, cube, tetrahedron, and prism, along with a hypercube of arbitrary dimensions.}

\end{abstract}

\begin{keyword}
Subdomain Integration, Embedded Interface Integration, Equivalent Polynomial, Polylogarithm
\MSC[2020] 65D32, 11G55, 65M60
\end{keyword}
\end{frontmatter}

\section{Introduction}

Partial differential equation (PDE) solvers are ubiquitous among many engineering and applied mathematics practitioners. Today, many PDE solvers employ discontinuous functions, especially in fluid dynamics problems. These methods use discontinuous functions to distinguish different subdomains and ensure no extrinsic contributions incur while utilizing an arbitrary discontinuity. A few particular extensions of the Finite Element Method (FEM) using discontinuous functions are CutFEM or Extended FEM (XFEM), generalized FEM (GFEM), and nonlocal FEM. In XFEM and GFEM, an enrichment function, e.g., the Heaviside function, is employed to distinguish different domains defined by a common interface, alleviating cumbersome remeshing techniques \cite{aragon2010generalized,moes1999finite}. The nonlocal FEM implements a kernel function, represented by the step function, that ensures nonlocal contributions are zero outside of some specified region \cite{aulisa2021nmpde, delia2020acta,friswell2007non}. Venturing outside of the FEM, an example of a method that also employs discontinuous functions is the Volume of Fluid (VOF) method. The VOF method uses the characteristic function to determine what portion of a cell is occupied by a fluid \cite{AULISA20072301,aulisa2003mixed, hirt1981volume}. From the methods above, one can see the critical role of discontinuous functions in many of today's PDE solvers, all of which benefit from an accurate and efficient way of dealing with the integration of a discontinuous function.  

Discontinuous function integration can be cast as integration over several disjoint subdomains involving continuous functions, i.e., the region over which the integration occurs can be broken up into multiple subdomains where only continuous functions are defined. However, the boundary defining the subdomains is rarely trivial, and traditional integration schemes are not practical. Even invoking the divergence theorem in such cases has proven to be intractable for even simple geometries and discontinuities, as seen in \cite{joulaian2016numerical,ventura2006elimination}, since integration over the subdomain boundaries must be performed. There has also been working devoted to moment fitting approaches, such as in \cite{joulaian2016numerical}, which also rely on the divergence theorem. Another approach, which depends on the convexity of the region of integration, is presented in \cite{mousavi2011numerical}. Although in general, it cannot be expected that the region of integration is convex. \jon{Exact quadrature rules exist only for triangle and tetrahedral geometries, with planar cuts, and for low degree polynomials, specifically for quadratic polynomials, \cite{holdych2008quadrature}. In there, the authors suggest they can extend their work to a higher degree polynomial integration. Still, the automation seems challenging since every time a higher degree polynomial is considered, one must recompute a new set of quadrature points and rules.}

The most common approach to discontinuous function integration is an adaptive algorithm, i.e., an algorithm that uses a grid refinement technique to capture the discontinuity better and produce a more accurate approximation to the integral. However, adaptive methods still require extensive information about the boundaries of the subdomains and typically lead to high computational costs, ultimately slowing down the numerical PDE scheme. There have been several recent developments that deal with the issue of discontinuous function integration, avoiding expensive adaptive methods. Among these are the use of equivalent polynomials \cite{abedian2019equivalent,abedian2013performance,ventura2015equivalent,ventura2006elimination}, more specifically, a polynomial that replaces the discontinuous function in the integrand and yields an equivalent integral. Equivalent polynomial methods allow for integration of continuous functions over an entire region without the difficulty of discontinuous functions and, for line/plane discontinuities, the high computation cost of adaptive quadrature methods.   

The equivalent polynomial method was first introduced by Ventura in \cite{ventura2006elimination}, where equivalent polynomials were found analytically for simple geometries and discontinuities \cite{gasca2000polynomial,ma1996generalized,neidinger2019multivariate}. The ideas presented in \cite{ventura2006elimination} were limited to lower order elements, e.g., linear triangles and tetrahedra, and bi-linear quadrilaterals, as a consequence of using the divergence theorem to analytically calculate the coefficients of the equivalent polynomial. The difficulty behind this method is introduced when analytical integration is applied to a generic discontinuity since integration must be carried out on two subdomains. This causes severe restrictions when the dimension increases, resulting in impractical discontinuity considerations, even when the discontinuity is a hyperplane. Moreover, the automation becomes difficult when the degree of the polynomial and/or the dimension increases.
\jon{In \cite{kees2022simple}, the authors implemented the work of Ventura up to 4$^{th}$ degree polynomials in 3D for the integration of fluid-structure Nitsche-type cutFEM coupling terms. However, the extension to higher degree polynomial integration can be numerically difficult.} 

The work in \cite{ventura2015equivalent}, by Ventura and Benvenuti, builds on the core idea presented in \cite{ventura2006elimination}, i.e., the idea of replacing a discontinuous function with an equivalent polynomial. The work's limitation in \cite{ventura2006elimination} is overcome by using a regularized Heaviside function, which approaches the Heaviside function in the limit, in place of the Heaviside function. This regularized Heaviside function is continuous and differentiable for any value of the regularization parameter $\rho.$ The regularized Heaviside function allows one to perform analytical integration over the entire domain and then take the limit of the resulting expression when deriving the equivalent polynomial coefficients. The extended work in \cite{ventura2015equivalent} creates a more robust method by eliminating the need for analytical integration over arbitrary subdomains created by the discontinuity. This method requires equality between the integral of the regularized Heaviside function multiplied by some monomial and the integral of the equivalent polynomial multiplied by the same monomial. The highest degree of the monomial and the dimension dictate the size of the linear system that needs to be solved to recover the coefficients of the equivalent polynomial. Since the equation for the discontinuity appears in the regularized Heaviside function, the equivalent polynomial coefficients will be dependent upon the discontinuity and the regularization parameter. Large values of the regularization parameter can then be taken to approximate the Heaviside function. Automation of this method relies on numerical libraries to calculate the polylogarithm function, which naturally arises from integrating the regularized Heaviside function.
Using the regularized Heaviside function, one is left with expressions that involve a linear combination of polylogarithm functions of various orders. The two sources of error arising from the use of equivalent polynomials, as mentioned in \cite{ventura2015equivalent}, are the numerical evaluation of the polylogarithm and round-off error introduced by large values of the regularization parameter $\rho.$      

In \cite{abedian2019equivalent} the concept of equivalent polynomials was extended to incorporate Legendre polynomials, which give rise to very beneficial properties. The main idea is to represent the equivalent and element shape polynomials with Legendre polynomials. The properties of Legendre polynomials are then utilized to allow for analytical integration over specified squares in 2-D or cubes in 3-D. Hence the error incurred from this method is produced by a spacetree refinement algorithm for complex discontinuities.  It is stated in \cite{abedian2019equivalent} that the analytical integration results are the same as those in \cite{ventura2015equivalent} for a line or plane discontinuity; however, the implementation of the equivalent Legendre polynomials for the specified discontinuities lacks the ease of algorithmic automation for the analytical integration. 

\eugenio{In this work, we overcome these difficulties and provide closed-form recursive algebraic formulas to exactly evaluate the subdomain integral for any degree polynomial on parent finite element shapes cut by a planar interface. The proposed method completes the equivalent polynomial technique in \cite{ventura2015equivalent}. Similar to their work, ours eliminates the need to integrate over a specified, often intractable, subdomain $\Omega_i$ by integrating on the whole domain and using the Heaviside function as a weight. Again the regularized Heaviside function is substituted by the polylogarithm function. Still, here, we take advantage of its derivative and limiting properties, yielding a formulation that ultimately eliminates the need for evaluating the polylogarithm. As a result, both sources of error introduced by the polylogarithm in \cite{ventura2015equivalent} are removed. By utilizing the derivative properties of the polylogarithm and the relationship the derivative of the Heaviside function shares with the Dirac delta distribution, we derive the same type of formulas to exactly evaluate the embedded interface integral for any degree polynomial. We provide the pseudo-codes for the subdomain and embedded interface integrals on hypercubes, triangles, tetrahedra, and prisms cut by planar surfaces. Much attention has been given to algorithms that avoid overflow in computer arithmetic, using, when needed, alternative formulas derived to eliminate round-off errors. The recursive nature of the algorithms allows for full automation. All expressions are algebraic and easy to implement.} 

The outline of this paper is as follows. In Section 2, we discuss the properties of the polylogarithm, which are implemented in the paper. In Section 3, the closed-form expressions for the different elements are derived. The pseudo-codes for the n-dimensional cube, triangle, tetrahedron, and wedge can also be found in this section. Lastly, in Section 4, we provide some useful notes on the practical implementation of the equivalent polynomials. First we deal with the ill-conditioned Gram mass matrix and then we show that our algorithm can be executed offline, while new quadrature rules can be computed online at very little cost and to any accuracy using simple interpolation. \eugenio{In the Conclusion section, we refer to a new result we obtained for curved cuts and for now only available  in the Ph.D. thesis of the second author \cite{loftin2022exact}}.

All the algorithms developed in this article are implemented in FEMuS \cite{aulisa2010femus},
an in-house open-source finite element C++ library built on top of PETSc \cite{balay2012petsc} and publicly available on GitHub.

\section{Preamble}
The polylogarithm, $\Li_s(z)$ where $s,z\in \mC$ with $|z|<1$, can be defined as $$\Li_s(z) = \sum_{k=1}^{\infty} \frac{z^k}{k^s},$$ or in integral form as $$\Li_s(z) = \frac{z}{\Gamma(s)} \int_{0}^{\infty} \frac{x^{s-1}}{e^x - z} \,dx, $$
by analytic continuation, where $\Gamma(s)$ is the gamma function. The integral representation of $\Li_s(z)$ is analytic for $z \in \mC \setminus [1,\infty)$ and $\Re(s)>0$ \cite{dingle1957fermi,truesdell1945function}. When the above integral is replaced with an appropriate complex contour integral we can consider $s \in \mZ^{-} \cup \{0\} $\cite{truesdell1945function}. For the purpose of this paper we will only consider polylogarithms of the form $\Li_s(w)$, where $s \in \{ -1, 0, 1, ... \}$ and $w \in \mR $. The polylogarithm can be defined in closed form for $s=1,0,-1,...\,.$ Specifically, $$ \Li_0(w) = \frac{w}{1-w}.  $$ All identities in this paper are used when $\Li_s(w)$ is well defined. 

Two useful properties used throughout this paper are $\Li_s(-e^w) = -F_{s-1}(w)$, where $F_{s-1}(w)$ is the Complete Fermi-Dirac integral, and 
\eugenio{
\begin{equation}
\frac{ d \Li_s(-e^{\mu})}{d\mu} = \Li_{s-1}(-e^{\mu}) \label{derivativeOfLi}   
\end{equation} \cite{dingle1957fermi,rhodes1950fermi}. 
Two cases are of particular importance: $s=0$ and $s=-1$. For $s=0$, the polylogarithm $\Li_0$ is used to represent the Heaviside function $\U$, and for $s=-1$, the polylogarithm $\Li_{-1}$ is used to represent the Dirac delta distribution $\delta$. Namely, for any smooth level set function $G(\bm x)$
\begin{equation} \label{heaviside_equality}
\U(G(\bm x)) = -\lim_{t\rightarrow \infty}
\Li_0(-\exp(G(\bm x)t)),   
\end{equation}
where $\U$ is the Heaviside function with half-maximum convention 
\begin{equation}\U(G(\bm x)) = \left\{
\begin{matrix} \label{heaviside_hmax}
1 &\mbox{for} &G(\bm x) >0\\
0.5 &\mbox{for}& G(\bm x) = 0\\
0 &\mbox{for} &G(\bm x) < 0
\end{matrix}\right. \;,\end{equation}
and for any differentiable function $f(\bm x)$
\begin{equation}\label{dirac_equality1} 
\int_D f(\bm x) \, \delta(G(\bm x))  \|\nabla G\| d\bm x =- \lim_{t\rightarrow \infty} t \int_D f(\bm x)\,   \Li_{-1}(-\exp(G(\bm x)t)) \|\nabla G\| d \bm x.
\end{equation}
Equality \eqref{dirac_equality1} is the weak convergence of $-t \Li_{-1}(-\exp( G(\bm x)t)) \|\nabla G\|$  to the the Dirac delta distribution $\delta(G(\bm x)) \| \nabla G\|$  \cite{lax2014functional,onural2006impulse}. 
The proof is quite technical and it is 
given in Appendix A.
For a hyperplane level set function $G(\bm x) = \bm n \cdot \bm x + d $, with unit normal $\bm n$, we have $\nabla G = \bm n$, $\|\nabla G\|=1$, and Eq.~\eqref{dirac_equality1} further simplifies to
\begin{equation}\label{dirac_equality} 
\int_D f(\bm x)\,\delta(\bm n \cdot \bm x+d) d\bm x =- \lim_{t\rightarrow \infty} t \int_D f(\bm x)\,   \Li_{-1}(-\exp((\bm x\cdot \bm n+d)t)) d \bm x.
\end{equation}
The results of this paper should be implemented with consideration given to the aforementioned properties.} 

\eugenio{
The following proposition, which applies to a general domain, relates subdomain and embedded interface integrals to the limits of certain domain integrals of the polylogarithm functions $\Li_0$ and $\Li_{-1}$. This result is obtained by first connecting the subdomain and interface integrals to corresponding domain integrals, using the Heaviside function $\U$ and the Dirac delta distribution $\delta$ as weights, and then by using equality \eqref{heaviside_equality} between $\U$ and the limit of $\Li_0$ and equality \eqref{dirac_equality} between the integral of $\delta$ and the limit of the integral of $\Li_{-1}$, respectively.}
\begin{figure}[!h]
    \centering
    \includegraphics[scale=0.75]{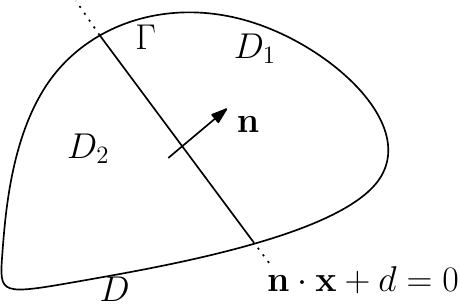}  \caption{Domain $D$ cut by the planar interface $\bm n \cdot \bm x +d = 0$.}
    \label{fig:Domain}
\end{figure}
\begin{proposition} \label{stepFunction}
\eugenio{With respect to Figure \ref{fig:Domain}, let the region $D$, with boundary $\partial D$, be cut by the plane $\bm n \cdot \bm x + d = 0$ into two subregions $D_1$ and $D_2$, with normal $\bm n$, pointing from $D_2$ to $D_1$. Let $\Gamma$ be the embedded interface between $D_2$ and $D_1$. Let $P_m(\bm x)$ be a polynomial of degree $m$ in $D$. Then, the following integral equalities hold}

\begin{align}
\int_{D_1} P_m(\bm x) d\bm x 
&=- \lim_{t \rightarrow \infty} \int_D P_m(\bm x) \Li_0(-\exp( (\bm n \cdot \bm x +d)t )) d \bm x \label{eq:volEquality}\\
&= \lim_{t \rightarrow \infty} \int_D P_m(\bm x) \left(1 + \Li_0(-\exp( -(\bm n \cdot \bm x +d)t ))\right) d \bm x . \label{eq:volEquality2}
\end{align}
Moreover, if $\Gamma$ is tangential to $\partial D$ at most on a set of measure zero and $\| \bm n\| = 1$, then
\begin{align} 
 \int_{\Gamma} P_m(\bm x) d \bm \mu 
&=-\lim_{t \rightarrow \infty} \int_D P_m(\bm x) \;t \Li_{-1}(-\exp( (\bm n \cdot \bm x  + d )t )) d \bm x\label{eq:surfEquality}\\
&=-\lim_{t \rightarrow \infty} \int_D P_m(\bm x) \;t \Li_{-1}(-\exp(-(\bm n \cdot \bm x  + d )t )) d \bm x\label{eq:surfEquality2}.
\end{align}
\end{proposition}

\begin{proof}
\eugenio{
Equality \eqref{eq:volEquality} follows from the integral equality
\begin{align*}
&\int_{D_1} P_m(\bm x) d \bm x =\int_D P_m(\bm x) \U(\bm n \cdot \bm x  + d) d\bm x,
\end{align*}
Eq.~\eqref{heaviside_equality}  and the dominated convergence theorem \cite{wang2018lecture}, i.e.,
\begin{align*}
\int_D P_m(\bm x) \U(\bm n \cdot \bm x  + d) d\bm x &=  \int_D P_m(\bm x) \left(-\lim_{t \rightarrow \infty} \Li_0(-\exp( (\bm n \cdot \bm x  + d) t))\right) d \bm x\\ 
 &= - \lim_{t \rightarrow \infty} \int_D P_m(\bm x) \Li_0(-\exp( (\bm n \cdot \bm x  + d) t)) d \bm x. 
\end{align*}
Eq. \eqref{eq:volEquality2} follows from Eq.~\eqref{eq:volEquality} and the following integral equality
\begin{align*}
  \int_{D_1} P_m(\bm x)d \bm x &=  \int_{D} P_m(\bm x) \left(1 - \U(-(\bm n \cdot \bm x  + d))\right) d \bm x.
\end{align*}
Equality \eqref{eq:surfEquality}  follows from the integral equality
\begin{align*}
  &\int_{\Gamma} P_m(\bm x) \bm d \mu =  \int_{D} P_m(\bm x) \delta(\bm n \cdot \bm x  + d) d \bm x,
  \end{align*}
and from Eq.~\eqref{dirac_equality}
\begin{align*}
& \int_D P_m(\bm x) \delta(\bm n \cdot \bm x  + d) d \bm x=- \lim_{t \rightarrow \infty} \int_D P_m(\bm x)\,  t \, \Li_{-1}(-\exp( (\bm n \cdot \bm x  + d)t )) d\bm x.
\end{align*}
Eq.~\eqref{eq:surfEquality2} follows from Eq.~\eqref{eq:surfEquality} and the following integral equality
\begin{align*}
 & \int_{D} P_m(\bm x) \delta(\bm n \cdot \bm x  + d) d \bm x =
 \int_{D} P_m(\bm x) \delta(-(\bm n \cdot \bm x  + d)) d \bm x.
\end{align*}
}
\end{proof}

\begin{remark}
In proving \eqref{eq:surfEquality}, we assumed that the interface $\Gamma$ is tangential to the boundary of $D$ only on a set of measure zero. Such distinction is needed, since otherwise \eugenio{the domain of the Dirac delta distribution}, centered on $\partial D$ and aligned with the normal direction, would be only half contained within $D$, thus contributing only for half to the interface integral. In all the applications we are going to consider next, $D$ will only be a convex domain with piece-wise flat boundaries. In doing so, $\Gamma$ is either completely tangential or never tangential to $\partial D$. 
This allows us to compute the interface integral also in the tangential case (the boundary integral) by doubling the value of the computed integral in \eqref{eq:surfEquality}.
\end{remark}

\eugenio{
In the remaining part of this section we will build the needed tools to exactly evaluate integrals as the ones in Eqs.~\eqref{eq:volEquality}-\eqref{eq:surfEquality2} for the one dimensional case, where the domain $D$ is the interval $[0,1]$ and level set function $G(x) = a x +d$. The resulting formula will hold for all integer $s\ge-1$.}

The following two propositions are a direct consequence of the properties found in \cite{dingle1957fermi,rhodes1950fermi,truesdell1945function}.

\begin{proposition}\label{limit}
For $s=0,1,2,\dots$, $$\LLi_s (a) := \lim_{t \rightarrow \infty} \frac{\Li_s (- \exp( a t) )}{t^s} =
\begin{cases}
-0.5 &\mbox{if } s=0 \mbox{ and } a = 0\\
-\dfrac{a^s}{s!} & \mbox{if } a > 0\\
0 & \mbox{otherwise}
\end{cases}
 .  $$
\end{proposition}
\begin{proof} 
\textbf{Case 1}: $s=a=0$. In this case we have $$\LLi_0 (0) =  \Li_0(-1) = \frac{-1}{1+1} = -\frac{1}{2}.$$

\textbf{Case 2}: By induction on $s$, with $a>0$.\\ 
For $s=1$ we have $$\LLi_1(a) = \lim_{t \rightarrow \infty} \frac{\Li_1 (- \exp( a t) )}{t} = \lim_{t \rightarrow \infty} \frac{-\ln{(1 + e^{at})}}{t} = -a.$$ 
\\
For $s=k-1$: assume $$ \lim_{t \rightarrow \infty} \frac{\Li_{k-1}(-e^{a t})}{t^{k-1}} = -\frac{a^{k-1}}{(k-1)!}.$$
\\
Then for $s=k$ we have $$ 
\lim_{t \rightarrow \infty} 
\frac{\Li_{k}(-e^{a t})}{t^{k}} =  \lim_{t \rightarrow \infty} 
\frac{a}{k} 
\frac{\Li_{k-1}(-e^{a t})}{t^{k-1}} 
= -\frac{a^{k}}{k!},
$$
\\
where we used \eugenio{the derivative property in Eq.~\eqref{derivativeOfLi}.}\\

\textbf{Case 3}: For any other case, i.e., $a \leq 0$, $s \ne 0$, the terms of the appropriate series expansion for the Fermi-Dirac integral vanish when the limit is taken inside the series, where the series converges uniformly\cite{dingle1957fermi}.

\end{proof}
\eugenio{In the following proposition we show that, for any integer $s\ge-1$, it is possible to derive close form expressions for the antiderivative of the polylogarithm function $\Li_s (-\exp(a x+d)t)$ multiplied by any monomial.}
\begin{proposition} \label{prop:defInt} Let $a \ne 0$, $m = 0,1,2, \dots$ and $s=-1,0,1,\dots$, then
$$ \int x^m \Li_s (-\exp(a x+d)t)\, dx  = \sum_{i=1}^{m+1} \frac{m! \, {(-1)}^{i-1}}{(m+1-i)!}  x^{m+1-i} \frac{\Li_{s+i} (-\exp (a x +d)t )}{(a t)^{i}} + C.$$

\end{proposition}

\begin{proof}
In the domain of interest, where the polylogarithm function converges uniformly, we use the identity
 $$ \Li_s(-e^{\mu}) = \int \Li_{s-1}(-e^{\mu}) d\mu+C, $$ 
 which follows from \eugenio{the derivative property in Eq.~\eqref{derivativeOfLi}.}
 For ease of notation we drop the constant in the proof. \\
For $m=0$, we get $$\int  \Li_s (-\exp( (a x +d )t )) dx = \frac{\Li_{s+1}(-\exp((a x + d)t))}{at}= \sum_{i=1}^{1} \frac{(-1)^{i+1}0!x^{0+1-i}\Li_{s+i}(-\exp((a x + d)t)) }{(0-i+1)!(ta)^i i^{s+1}}.   $$ 
For $m=k$ assume 
$$ \int x^k \Li_{s+1} (-\exp( (a x +d )t )) dx  = \sum_{i=1}^{k+1} \frac{k! \ {(-1)}^{i+1}}{(k+1-i)!}  x^{k+1-i} \frac{\Li_{s+1+i} (-\exp( (a x +d )t )}{(a t)^i}.$$ 
\
Then for $m=k+1$, from integration by parts, we have
\begin{align*}
\int& x^{k+1} \Li_s (-\exp( (a x +d )t )) dx  \\
&= x^{k+1} \int \Li_s (-\exp( (a x +d )t )) dx - (k+1) \int x^k\frac{\Li_{s+1} (-\exp( (a x +d )t   )}{at}dx\\  
&=  \frac{x^{k+1}}{at} \Li_{s+1} (-\exp( (a x +d )t )) - (k+1)\sum_{i=1}^{k+1} \frac{k! \ {(-1)}^{i+1}}{(k+1-i)!}  x^{k+1-i} \frac{\Li_{s+i+1} (-\exp( (a x +d )t )}{(a t)^{i+1}}  \\ &=  \sum_{i=1}^{k+2} \frac{(k+1)! \ {(-1)}^{i+1}}{(k+2-i)!}  x^{k+2-i} \frac{\Li_{s+i} (-\exp( (a x +d )t )}{(a t)^{i}}. \\  
\end{align*} 
This completes the proof.
\end{proof}

\eugenio{
Integrals, as the ones in Proposition \ref{stepFunction}, require the evaluation of the limit for $t\rightarrow\infty$ of the definite integral
of polynomials multiplied by polylogarithm functions. Using the results in Propositions \ref{limit} and \ref{prop:defInt}
we are now capable of evaluating such limits for one-dimensional integrals. Later, in Section \ref{section3}, we will extend these results to limits of integrals in higher dimensions.}

\eugenio{
\begin{proposition}
For , $m = 0,1,2, \dots$ and $s=-1,0,1,\dots$, the limit of the following definite integral is given by
\begin{align}
I_1 &= -\lim_{t \rightarrow \infty} \frac{1}{t^s}\int_0^1 x^m \Li_s(-\exp( (a x +d )t ) \,dx \nonumber \\
&=\sum_{i=1}^{m+1} \frac{m!}{(m+1-i)!} \frac{1}{(-a)^i} 
\LLi_{s+i}(a+d) 
- \frac{m!}{(-a)^{m+1}} \LLi_{s+m+1} (d) .\label{int2}
\end{align}
\end{proposition}
}
\begin{proof}
The proof follows directly from combining Propositions \ref{limit} and \ref{prop:defInt}.
\end{proof}

\begin{definition}
For $a\ne 0$, let
\begin{equation}
 I_2 = \sum_{i=0}^{s} \frac{(-a)^{s-i}  (a+d)^i}{i!(m+1+s-i)!}.\label{int2a}
\end{equation}
\end{definition}

\begin{remark} \label{remarkOverflow1} \eugenio{In the following proposition, we will show that
Eq.~\eqref{int2} is equivalent to Eq.~\eqref{int2a} for $s\ge 0$ and positive arguments of the $\LLi$ functions.} In  computer arithmetic Eq.~\eqref{int2} suffers from overflow for $d\gg|a|>0$ and $a\rightarrow0$, because of the presence of the $\frac{1}{a^i}$ terms in the sums. Proposition \ref{equivalentFormula} will show all these terms actually simplify after expanding the definition of $\LLi$ for positive argument.
\end{remark}

\begin{remark}\label{remarkOverflow2}
In  Eq.~\eqref{int2}, for $d\le0$ or $ a+d \le 0$, and $a\rightarrow0$, either the arguments of the polylogarithm functions are non positive, or, if positive, they are of the same order of $a$. In the first case  the contribution of their limits is zero. In the second case using the definition of $\LLi$ with positive argument one would get
$$ \frac{\LLi_{s+i}(O(a))}{a^i} =  \frac{O(a)^{s+i}}{a^i}= O(a)^s$$ for all $i$. Thus, all terms in Eq.~\eqref{int2} would have comparable size and, since $a$ does not appear in the denominator, it no longer contributes to overflows for $a\rightarrow0$.
\end{remark}

\begin{proposition} \label{equivalentFormula} $\,$\\
For $s\ge0$, $a\ne0$, $d>0$, and $a+d>0$ $$I_1 = I_2.$$
\end{proposition}

\begin{proof}
First note that the conditions $a\ne0$, $d>0$ and $a+d>0$ are equivalent to $\dfrac{-d}{a} \not\in [0,\, 1]$.

In proving the proposition one simply needs to apply integration by parts and utilize Proposition \ref{limit}. For a fixed $s\ge0$ and $a\ne0$ we have
\begin{align*}
-\lim_{t \rightarrow \infty}& \frac{1}{t^s}\int_0^1 x^m \Li_s(-\exp( (a x +d )t ) \,dx  \nonumber \\ &= -\lim_{t \rightarrow \infty} \frac{1}{t^s}\left( \frac{x^{m+1}}{m+1} \Li_s(-\exp( (a x +d )t )\Big|_0^1 -\frac{at}{m+1}\int_0^1 x^{m+1} \Li_{s-1}(-\exp( (a x +d )t)) \, dx\right)\nonumber \\ 
&=\frac{(a+d)^s}{s!(m+1)} - \frac{a}{m+1}\, \left(-\lim_{t \rightarrow \infty} \frac{1}{t^{s-1}}\int_0^1 x^{m+1} \Li_{s-1}(-\exp( (a x +d )t)) \, dx\right) \nonumber\\
&=\frac{(a+d)^{s}m!}{s!(m+1)!} - \frac{a(a+d)^{s-1}m!}{(s-1)!(m+2)!}  \nonumber \\
&\hspace{3cm}-\frac{a^2}{(m+1)(m+2)}\left(-\lim_{t \rightarrow \infty} \frac{1}{t^{s-2}}\int_0^1 x^{m+2} \Li_{s-2}(-\exp( (a x +d )t)) \, dx\right) \nonumber \\
&=\quad \dots \nonumber \\
&=\sum_{i=0}^{s} \frac{m!(-a)^{s-i}  (a+d)^i}{i!(m+1+s-i)!} -\frac{(-a)^{s+1}m!}{(m+2+s)}\left(-\lim_{t \rightarrow \infty} \int_0^1 x^{m+2+s}\, t\, \Li_{-1}(-\exp( (a x +d )t)) \, dx\right)  \nonumber\\
&=\sum_{i=0}^{s} \frac{m!(-a)^{s-i}  (a+d)^i}{i!(m+1+s-i)!},
\end{align*}
where we have used the weak convergence of $ t  \Li_{-1}(-\exp( (ax  + d)t ))$  to the the Dirac distribution $\delta$ with $\frac{-d}{a} \not\in [0,\, 1]$.
\end{proof}

\section{Polynomial basis integration} \label{section3}
In this section, we provide closed-form algebraic expressions for the integration of selected polynomial bases for several standard FEM shapes. For a particular element, the polynomial basis implemented results in a monomial integrand, after an appropriate transformation. For each element, several integration strategies are provided in order to avoid overflow in computer arithmetic. Each integration is then summarized into a detailed algorithm.
\subsection{LSI: Line Segment Integration on $[0,1]$, with $a \ne 0$}
For a fixed $s=-1,0,1,\dots$, we want to evaluate integrals in the form
$$ \LSI_s^m(a,d) = -\lim_{t\rightarrow \infty}  \frac{1}{t^s} \int_{0}^{1} x^m \Li_s(-\exp( (a x + d )t ) dx.$$
From Eq.~\eqref{int2}
\begin{align} \LSI_s^m(a,d) &:= \sum_{i=1}^{m+1} \frac{m!}{(m+1-i)!} \frac{1}{(-a)^i} 
\LLi_{s+i}(a+d) 
- \frac{m!}{(-a)^{m+1}} \LLi_{s+m+1} (d) . \label{lineFormula}
\end{align}
For all $a$, $d\in\mathbb{R}$, with $|a|>0$, we have the subdomain integral
\begin{equation}\label{LSI0}
\LSI_0^m(a,d) = \int_{0}^{1} x^m \U(ax+d) dx,
\end{equation}
and, for $a^2=1$, the interface integral
\begin{equation} \label{LSIm1a}
\LSI_{-1}^m(a,d) = \int_{0}^{1} x^m \delta(ax+d) dx.
\end{equation}
For $s=-1$ and $|a|>0$ we also have the explicit point evaluation formula
\begin{equation} \label{LSIm1}
\LSI_{-1}^m(a,d) = \begin{cases}
\vspace{3pt}
\displaystyle \frac{1}{|a|}\left(-\frac{d}{a}\right)^m  &  \ \mbox{if } 0<-\dfrac{d}{a} < 1\\ \vspace{3pt}
\displaystyle \dfrac{1}{2|a|} \left(-\frac{d}{a}\right)^m  &  \ \mbox{if }\, -\dfrac{d}{a} =0\, \mbox{ or }\,-\dfrac{d}{a} =1 \\ \vspace{3pt}
0 &  \ \mbox{elsewhere}
\end{cases},   
\end{equation}
with the assumption that $0^0=1$. That is the case for $m=0$ and $d=0$.
Although equivalent to Eq.~\eqref{lineFormula}, for $s=-1$, Eq.~\eqref{LSIm1} is generally faster to compute and does not suffer from overflow in computer arithmetic.
Also, for $s\ge0$, $d>0$, and $a+d>0$ we replace Eq.~\eqref{lineFormula} with the equivalent Eq.~\eqref{int2a} to avoid overflow. 
The pseudo-code for the line segment integration is given in Algorithm \ref{algLSI}.

\begin{remark} \label{rmk:bd1}
The formula for $\LSI_{-1}^m(a,d)$ halves the value of the interface integral if the point $-d/a$ is one of the two boundary points. This happens because half of the domain of the Dirac distribution falls outside the line segment, thus it does not contribute to the integral value. If this is not the desired behavior, and the boundary integral should account for the whole value, the definition of $\LSI_{-1}^m(a,d)$ should be replaced by
\begin{equation} \label{LSIm1b}
\LSI_{-1}^m(a,d) = \begin{cases}
\vspace{3pt}
\displaystyle \frac{1}{|a|}\left(-\frac{d}{a}\right)^m  &  \ \mbox{if } 0\le-\dfrac{d}{a} \le 1\\ 
0 &  \ \mbox{elsewhere}
\end{cases},   
\end{equation}
again with the assumption that $0^0=1$.
\end{remark}
\begin{algorithm} 
\caption{Pseudo-code for integration on the line segment $[0,1]$ with $a\ne0$ and $s=-1,0,1,\dots$. For $s=-1$ and $|a|=1$ it corresponds to the interface integral.
For $s=0$  it corresponds to the subdomain integral.}
\label{algLSI}
\begin{algorithmic}[1] 
  \Function {Line\_Segment\_Integration}{$a$, $d$, $m$, $s$}
  \If{$s=-1$}
    \State{\Return $\; \LSI_{-1}^m(a,d) $ from Eq.~\eqref{LSIm1}} 
  \Else   
     \If {$d \le 0$ or $a+d \le 0$}
     \State{\Return $\qquad \mathlarger{\sum}_{i=1}^{m+1} \dfrac{m!}{(m+1-i)!} \dfrac{1}{(-a)^i} 
\LLi_{s+i}(a+d) 
- \dfrac{m!}{(-a)^{m+1}} \LLi_{s+m+1} (d) $}    
     \Else
     \State{\Return $ \qquad \mathlarger{\sum}_{i=0}^{s} \dfrac{m!(-a)^{s-i}  (a+d)^i}{i!(m+1+s-i)!}$}
      \EndIf
    \EndIf
  \EndFunction
  \end{algorithmic}
\end{algorithm}

\subsection{SQI: Square Integration on $[0,1]^2$, with $a^2+b^2>0$}
Fix $s=-1,0,1,\dots$, we want to evaluate integrals in the form
$$ \SQI_s^{mn} (a,b,d) = -\lim_{t\rightarrow \infty}  \frac{1}{t^s} \int_{0}^{1} \int_{0}^{1} x^m y^n \Li_s(-\exp( (a x +b y + d )t ) dy\,dx.$$
We will first consider the case when the interface $\Gamma$ is parallel to either the square sides, and then all the remaining cases.

If $a=0$ the iterated integral can be split in the product of two integrals
$$ \SQI_s^{mn} (0,b,d) = \int_{0}^{1} x^m dx \left(-\lim_{t\rightarrow \infty}  \frac{1}{t^s}  \int_{0}^{1} y^n \Li_s(-\exp( (b y + d )t ) dy \right) = \frac{1}{m+1} \LSI_s^n(b,d).$$

Similarly, if $b=0$
$$ \SQI_s^{mn} (a,0,d)  = \frac{1}{n+1} \LSI_s^m(a,d).$$

If both $a$ and $b$ are different from zero, after the integration of the inner integral
we get
\begin{align}
\SQI_s^{mn}  (a,b,d) &= 
-\lim_{t\rightarrow \infty}  \int_{0}^{1} x^m \left(
-\sum_{j=1}^{n+1} \frac{n!} {(n+1-j)!} \, \frac{1}{(-b)^j} \frac{\Li_{s+j} (-\exp( (a x+b +d )t )}{ t^{s+j}}  \right.\nonumber \\
&\hspace{1cm}\left.+ \frac{n!}{(-b)^{n+1}} \frac{\Li_{s+n+1} (-\exp( (ax + d )t )}{ t^{s+n+1}} \right) dx \nonumber\\
&= -\sum_{j=1}^{n+1} \frac{n!} {(n+1-j)!} \, \frac{1}{(-b)^j}    \LSI_{s+j}^m(a,b+d) + \frac{n!}{(-b)^{n+1}} \LSI_{s+n+1}^m(a,d) , \label{iterative1}
\end{align}
Then, for all $a$, $b$, $d\in\mathbb{R}$ , such that $a^2+b^2>0$, we have the subdomain integral
$$\SQI_0^{mn}(a,b,d) = \int_{0}^{1} \int_{0}^{1} x^m y^n \U(ax+by+d) \,dy\,dx,$$
and, for $a^2+b^2 = 1$, the interface integral
$$\SQI_{-1}^{mn}(a,b,d) = \int_{0}^{1} \int_{0}^{1} x^m y^n \delta(ax+by+d)\, dy\,dx.$$
These formulas are general and versatile: they work regardless of where the line $a x + b y + d = 0$ intersects the square domain, and the orientation of the Heaviside function follows the orientation of the normal $\left<a,b\right>$.

\begin{remark} \label{rmk:bd2}
In the special cases $\SQI_{-1}^{mn}(a,0,d)$ (or $\SQI_{-1}^{mn}(0,b,d$)), with $-d/a= 0 \mbox{ or }1 $ (or $-b/d=0 \mbox{ or }1$), the corresponding line $ax+d=0$ (or $by+d=0$) overlaps with one of the sides of the square. Depending on which definition is used for $\LSI_{-1}^{m},$ either Eq.~\eqref{LSIm1} or Eq.~\eqref{LSIm1b}, one is left with half the boundary integral or the entire boundary integral, respectively, over the specified side of the square. This is also the case for the cube and the hypercube we are going to consider next. 
\end{remark}

\subsection{CBI: Cube integration on $[0,1]^3$, with $a^2+b^2+c^2>0$.}

Fix $s=-1,0,1,\dots$, we want to evaluate integrals in the form
$$ \CBI_s^{mno} (a,b,c,d) = -\lim_{t\rightarrow \infty}  \frac{1}{t^s} \int_{0}^{1} \int_{0}^{1} \int_{0}^{1} x^m y^n z^o \Li_s(-\exp( (a x +b y +c z+ d )t ) dz\,dy\,dx.$$
For $c=0$ the above integral reduces to the square case, i.e.,
\begin{align*}
& \CBI_s^{mno}(a,b,0,d) = \frac{1}{o+1} \SQI_s^{mn}(a,b,d).
\end{align*}
If $c\ne0$ and both $a=0$ and $b=0$
\begin{align}
\CBI&_s^{mno} (0,0,c,d) = \frac{1}{m+1} \frac{1}{n+1}\LSI_s^o(c,d).
\end{align}
For all other cases, after integrating in $z$ we get
\begin{align}
\CBI_s^{mno} (a,b,c,d)& 
= -\lim_{t\rightarrow \infty}  \int_{0}^{1}\int_{0}^{1} x^m y^n \left(-
\sum_{k=1}^{o+1} \frac{o!} {(o+1-k)!} \, \frac{1}{(-c)^k} 
 \frac{\Li_{s+k} (-\exp( (a x+by+c +d )t )}{ t^{s+k}} \right. \nonumber  \\
& \hspace{1cm}\left. +\frac{o!}{(-c)^{o+1}} \frac{\Li_{s+o+1} (-\exp( (ax +by + d )t )}{ t^{s+k}} \right) \, dy \, dx \nonumber\\
&= -\sum_{k=1}^{o+1} \frac{o!} {(o+1-k)!} \, \frac{1}{(-c)^k} \SQI_{s+k}^{mn}(a,b,c+d) + \frac{o!}{(-c)^{o+1}}  \SQI_{s+o+1}^{mn}(a,b,d). \label{iterative2}
\end{align}
The cases $a=0$ or $b=0$ are handled by the square integrals as described in the previous section. 

Then, for all $a$, $b$, $c$, $d\in\mathbb{R}$, with $a^2+b^2+c^2>0$, we have the subdomain integral
$$\CBI_0^{mno}(a,b,c,d) = \int_{0}^{1} \int_{0}^{1} \int_{0}^{1} x^m y^n z^o \U(ax+by+cz+d) \,dz\,dy\,dx,$$
and, for $a^2+b^2+c^2 = 1$, the interface integral
$$\CBI_{-1}^{mno}(a,b,c,d) = \int_{0}^{1} \int_{0}^{1} \int_{0}^{1} x^m y^n z^o \delta(ax+by+cz+d) \, dz\,dy\,dx.$$

It is remarkable how such simple formulas can handle all possible intersections between the cube and the plane. Moreover, they can be easily extended to evaluate corresponding integrals on hypercubes cut by hyperplanes for any dimension.

\subsection{HCI: Hypercube Integration on $[0,1]^{dim}$, with $\bm n =\langle a_1,a_2,\dots,a_{dim} \rangle, \|\bm n\|>0$ and $\bm m =\langle m_1,m_2,\dots,m_{dim} \rangle$.}
We are seeking integrals in the form 
$$ \HCI_{s,dim}^{\bm m} (\bm n,d) = -\lim_{t\rightarrow \infty}  \frac{1}{t^s} \int_{[0,1]^{dim}} \prod_{i=1}^{dim} x_i^{m_i} \Li_s(-\exp( (\bm n \cdot \bm x + d )t ) d \bm x,$$
where we assume $|a_i| \le |a_{i+1}|$. However, if this is not the case, one can perform a reordering of the normal coefficients due to the symmetry of the domain and the integrand. Define $dim_0 \in \mints_0$ with $dim_0 \leq dim$ such that $dim_0$ is an upper bound for the indices corresponding to all the $a_i =  0$  $\forall i < dim_0$. Define $dim' := dim - dim_0$, $\bm m' := \langle m_{dim_0+1},\dots, m_{dim}  \rangle$ and $\bm n' :=\langle a_{dim_0+1},\dots, a_{dim}  \rangle  $.
Then
$$ \HCI_{s,dim}^{\bm m} (\bm n,d) = \prod_{i=1}^{dim_0} \frac{1}{1 + m_i} \HCI_{s,dim'}^{\bm m'} \big( \bm n' ,d \big).$$
Then, dropping the $'$ superscript, the problem reduces to evaluating integrals in the form
$$ {\HCIA}_{s,dim}^{\bm m} (\bm n,d) = -\lim_{t\rightarrow \infty}  \frac{1}{t^s} \int_{HC_{dim}} \prod_{i=1}^{dim} x_i^{m_i} \Li_s(-\exp( (\bm n \cdot \bm x + d )t ) d \bm x,$$
with $|a_i| \le |a_{i+1}|$ and $a_1 \ne 0$. 
Following the same integration strategy used for the square and the cube, with $m = m_{dim}$ and $a=a_{dim}$, we obtain the following recursive formula
\begin{align}
{\HCIB}_{s,dim}^{\bm m} (\bm n,d) &= -\sum_{i=1}^{m+1} \frac{m!} {(m+1-i)!} \, \frac{1}{(-a)^i} 
{\HCIA}_{s+i,dim-1}^{\bm m^-} (\bm n^-, a + d)  \nonumber \\
&\hspace{1cm}+\frac{m!}{(-a)^{m+1}}  {\HCIA}_{s+m+1 , dim-1}^{\bm m^-}(\bm n^-,d),
\end{align}
where $\bm m^- = \langle m_{1},\dots, m_{dim-1}  \rangle$ and $\bm n^{-} =\langle a_{1},\dots, a_{dim-1}  \rangle  $.
This formula is recursively applied until dimension 1, where the the line segment integration formula, $\LSI$, is used.
At each level of integration two contributions occur, one that involves a sum and a single term. The most expensive terms to compute 
are the ones involving a summation, 
with each one of them requiring the computation of
$$\LLi_{s+k}\left(\sum_{i=1}^{dim} a_i +d\right),$$ for some $k\ge dim$. It is then desirable to have $$ \sum_i^{dim} a_i +d < 0,$$ so that all the $\LLi$ contributions vanish. From Proposition \ref{stepFunction},  changing the sign of the normal without any contribution is only allowed for $s=-1$, hence
$$ {\HCIB}_{-1,dim}^{\bm m} (\bm n,d) ={\HCIB}_{-1,dim}^{\bm m} (-\bm n,-d).$$
Similarly to Remarks \ref{remarkOverflow1} and \ref{remarkOverflow2}, the 
${\HCIB}$ formula also suffers from overflow in computer arithmetic when
$$\sum_{i=1}^{dim} a_i +d>>|a_{dim}|.$$
To overcome these difficulties, for this case only, we introduce  the alternative formula
\begin{align}
{\HCIC}_{s,dim}^{\bm m} (\bm n,d) &= \sum_{i=0}^{s} \frac{m!} {(m+1+i)!} \, {(-a)^i} 
{\HCIA}_{s-i,dim-1}^{\bm m^-} (\bm n^-, a + d)  \nonumber \\
&\hspace{1cm} + \frac{m!}{(m+s+1)!} (-a)^{s+1}  {\HCIA}_{-1 , dim}^{\bm m^*}(\bm n,d),
\end{align}
where $\bm m^* = \langle m_1,m_2,\dots,m_{dim-1},m + s +1 \rangle$.
This formula is obtained from \eugenio{the derivative identity in Eq.~\eqref{derivativeOfLi} and by recursive integration by parts}, increasing the monomial power and reducing the polylogarithm order $s$ until it reaches $-1$.
More specifically, ${\HCIC}$ follows the idea in Proposition \ref{equivalentFormula}, where an equivalent closed-form expression is given in which $a$ does not appear in the denominator. Note that in ${\HCIC}$ 
$$ {\HCIA}_{-1 , dim}^{\bm m^*}(\bm n,d)  = {\HCIA}_{-1 , dim}^{\bm m^*}(-\bm n,-d),$$
which not only permits choosing the optimal sign for the normal $\bm n$, but also
satisfies 
$$\sum_{i=1}^{dim} (-a_i) -d<<-|a_{dim}|.$$

The pseudo-code for general dimension $dim \ge 1$ is given in Algorithms \ref{alg1} and \ref{alg2}. In Algorithm \ref{alg1}, the contributions of each component with a zero coefficient $a_i$ are handled first. Algorithm \ref{alg2} is then called to compute the contributions from all the remaining components. The recursive nature of the algorithm follows from the patterns developed in the $\HCIA$, $\HCIB$, and $\HCIC$ formulas. 

Note that Algorithm~\ref{alg1} also handles the case $\bm n =\bm 0$. Although this case was excluded here, it will be needed later when integrating on the prism.

\begin{algorithm} 
\caption{Pseudo-code for the integration on the hypercube $[0,1]^{dim}$ cut by the hyperplane $\bm n \cdot \bm x + d=0$ with  $\bm n = \left<a_1,a_2,\dots,a_{dim}\right>$, $\bm m = \left<m_1,m_2,\dots,m_{dim}\right>$ and $s=-1,0,1,\dots$. For $s=-1$ and $\|\bm n\|=1$ it corresponds to the interface integral. 
For $s=0$ it corresponds to the subdomain integral.}
\label{alg1}
\begin{algorithmic}[1] 
  \Function {Hypercube\_Integration}{$dim$, $\bm n$, d, $\bm m$, $s$}
    \State{$\HCI = 1$}
    \For{ $i=1,\dots,dim$}
     \If {$a_i = 0$}
       \State{$\HCI \timeseq \dfrac{1}{m_i+1}$} 
       \State{Remove the $i$-th component of $\bm n$ and $\bm m$}
       \State{$dim = dim -1$}
       \State{$i = i - 1$}
     \EndIf
    \EndFor
    \If {$dim>0$}
    \State{Sort $\bm n$, and accordingly $\bm m$, from the smallest to the largest coefficient in magnitude}
    \State{\Return $\HCI * \text{\sc{Hypercube\_Integration\_A}}(dim,\, \bm n,\, d,\, \bm m,\, s)$}
    \Else \State {\Return $-\HCI * \LLi_s(d)$ }
    \EndIf
  \EndFunction
  \end{algorithmic}
\end{algorithm}
\begin{algorithm} 
\caption{Pseudo-code for the integration on the hypercube $[0,1]^{dim}$ cut by the hyperplane $\bm n \cdot \bm x + d=0$ with  $\bm n = \left<a_1,a_2,\dots,a_{dim}\right>$, $a_1 \ne 0$ and $|a_i|\le|a_{i+1}|$ for all $i=1,...,dim-1$, $\bm m = \left<m_1,m_2,\dots,m_{dim}\right>$, and $s=-1,0,1,\dots\;$. For $s=-1$ and $\|\bm n\|=1$ it corresponds to the interface integral. 
For $s=0$ it corresponds to the subdomain integral.}

\label{alg2}
\begin{algorithmic}[1] 
  \Function {Hypercube\_Integration\_A}{$dim$, $\bm n$, $d$, $\bm m$, $s$}
    
    \If{$dim = 1$} 
        \State{\Return {\sc{Line\_Segment\_Integration}}$(a_1\, d,\, m_1,\, s)$}
    \EndIf
  
    \State{$sum = \sum_{i=1}^{dim} a_i +d$}
    \If{$s = -1$}
        \If {$sum \le 0$} 
            \State {\Return $\text{\sc{Hypercube\_Integration\_B}} (dim,\, \bm n,\, d,\, \bm m,\, -1)$}
        \Else
            \State{\Return $\text{\sc{Hypercube\_Integration\_B}} (dim,\, -\bm n,\, -d,\, \bm m,\, -1)$}
        \EndIf
    \Else
        \If{ $sum \le |a_{dim}|$}
            \State{\Return $\text{\sc{Hypercube\_Integration\_B}} (dim,\, \bm n,\, d,\, \bm m,\, s)$}
        \Else
            \State{\Return $\text{\sc{Hypercube\_Integration\_C}} (dim,\, \bm n,\, d,\, \bm m,\, s)$}
        \EndIf
  \EndIf
  
  \EndFunction
  \end{algorithmic}
  \vspace{0.25cm}
\begin{algorithmic}[1] 
  \Function {Hypercube\_Integration\_B}{$dim$, $\bm n$, $d$, $\bm m$, $s$}
        \State {$m = m_{dim}; \quad a = a_{dim}$}
        \State{Remove the last component of $\bm n$ and  $\bm m$}
        \State {\Return 
        \begin{align*}\hspace{.75cm}
         &-\sum_{i=1}^{m+1} \frac{m!} {(m+1-i)!} \, \frac{1}{(-a)^i} 
         \text{\sc{Hypercube\_Integration\_A}}(dim-1,\, \bm n,\, a+d,\, \bm m,\, s + i)  \nonumber \\
         &\hspace{1cm}+\frac{m!}{(-a)^{m+1}}  \text{\sc{Hypercube\_Integration\_A}}(dim-1,\, \bm n,\, d,\, \bm m,\, s + m +1)
        \end{align*}}
  \EndFunction
  \end{algorithmic}
  \vspace{0.25cm}
  \begin{algorithmic}[1] 
  \Function {Hypercube\_Integration\_C}{$dim$, $\bm n$, $d$, $\bm m$, $s$}
        \State {$m = m_{dim}; \quad a = a_{dim}; \quad m_{dim} = m_{dim} + s + 1$}
        \State{$$\HCI = \dfrac{m!}{(m+s+1)!} (-a)^{s+1} \text{\sc{Hypercube\_Integration\_A}}(dim, \bm n,\, d,\, \bm m,\, -1)$$}
        \State{Remove the last component of $\bm n$ and $\bm m$}
        \State{\begin{align*}\hspace{.75cm}
          \HCI \pluseq \sum_{i=0}^{s}
          \dfrac{m!}{(m + i + 1)!} \, (-a)^{i}
            \text{\sc{Hypercube\_Integration\_A}}(dim-1, \bm n,\, a+d,\, \bm m,\, s - i) 
        \end{align*}}
  \State {\Return HCI }
  \EndFunction
  \end{algorithmic}
  
\end{algorithm}

\subsection{TRI: Triangle Integration 
on 
$\left\{ \protect\begin{array}{@{}l@{}} 0\le x\le 1 \\ 0\le y\le 1-x
         \protect\end{array} \right.,$
with $a^2+b^2>0$}

To ease the computation we choose a non-standard polynomial basis, namely $(1-x)^m y^n$. We then seek integrals
in the form
\begin{equation}
\TRI_s^{m,n} (a,b,d)= -\lim_{t \rightarrow \infty}\frac{1}{t^s}\iint_{\mbox{\sc{Tri}}} (1-x)^m y^n \Li_s(-\exp( (a x +b y +d )t )\, dA.    \label{TriOrigin}
\end{equation}
Then, for all $a$, $b$, $d\in\mathbb{R}$ such that $a^2+b^2 >0$, the subdomain integral is given by
$$\TRI_0^{mn}(a,b,d) = \iint_{\mbox{\sc{Tri}}} (1- x)^m y^n \U(ax+by+d) \, dA,$$
and, for $a^2+b^2 = 1$, the interface integral is given by $$\TRI_{-1}^{mn}(a,b,d) = \iint_{\mbox{\sc{Tri}}} (1-x)^m y^n \delta(ax+by+d) \, dA.$$

In Eq.\eqref{TriOrigin}, changing variables and renaming constants as follows
$$ x' = 1-x, \, y = y, \quad a'=-a, \,b'= b,\, d' = d + a,$$
yields
$$ -\lim_{t \rightarrow \infty}\frac{1}{t^s}\int_{0}^{1} \int_0^{x'} {x'}^{\,m} {y'}^{\,n} \Li_s(-\exp( (a' x' +b' y' +d' )t )) dy'\,dx'.$$
Dropping the $'$ superscript, for a fixed $s=-1,0,1,\dots$, the problem reduces to evaluating integrals in the form 
\begin{equation}
\TRIA_s^{mn}(a,b,d)=-\lim_{t \rightarrow \infty}\frac{1}{t^s}\int_{0}^{1} \int_0^x x^m y^n \Li_s(-\exp( (a x +b y +d )t ) dy\,dx. \label{TriAfterTransformation}    
\end{equation}

First, we will consider the three separate cases where the interface $\Gamma$ is parallel to one of the triangle edges.

For $b=0$
\begin{align}
\TRIB_s^{mn}&(a,0,d)=-\lim_{t \rightarrow \infty}\frac{1}{t^s}\int_{0}^{1} \int_0^x x^m y^n \Li_s(-\exp( (a x +d )t ) \,dy\,dx  \nonumber \\
&= -\lim_{t \rightarrow \infty}\frac{1}{t^s} \int_0^1 \frac{x^{m+n+1}}{n+1}  \Li_s(-\exp( (a x +d )t ) dx \nonumber \\
&= \frac{\LSI_s^{m+n+1}(a,d)}{n+1}.
\label{verticalCut}
\end{align}

For $a=0$ 
\begin{align}
\TRIB_s^{mn}&(0,b,d)=-\lim_{t \rightarrow \infty}
\frac{1}{t^s}\int_{0}^{1} \int_y^1 x^m y^n \Li_s(-\exp( (b y +d )t ) \,dx\,dy \nonumber \\
&= -\lim_{t \rightarrow \infty}\frac{1}{t^s} \frac{1}{m+1} \int_0^1 (y^n - y^{n+m+1}) \Li_s(-\exp( (b y +d )t ) dy \nonumber \\ 
&=\frac{\LSI_s^{n}(b,d)-\LSI_s^{m+n+1}(b,d)}{m+1}.
\label{horizontalCut}
\end{align}

For $ a + b = 0 $
\begin{align}
\TRIB_s^{mn}&(a,-a,d)=-\lim_{t \rightarrow \infty}\frac{1}{t^s}\int_{0}^{1} \int_0^x x^m y^n \Li_s(-\exp( (a x -a y +d )t ) \,dy\,dx \nonumber  \\
&= -\lim_{t \rightarrow \infty}\frac{1}{t^s} \int_0^1  n!
\left(\sum_{i=1}^{n+1} \frac{(-1)^{i-1}}{(n+1-i)!} x^{m+n+1-i} \frac{\Li_{s+i}(-\exp( (a x -a x+ d )t )}{(-at)^i} \right.  \nonumber \\
&\qquad\qquad\qquad\qquad\qquad\qquad\quad\left.- 
(-1)^{n} x^m \frac{\Li_{s+n+1}(-\exp(a x + d )t )}{(-at)^{n+1}}
\right) dx  \nonumber \\
&= n! \left(\sum_{i=1}^{n+1}  \left(\frac{1}{a}\right)^i \frac{ \LLi_{s+i}(d)}{(n+1-i)!} \int_0^1 x^{m+n+1-i} dx + \left(\frac{1}{a}\right)^{n+1}\LSI_{s+n+1}^{m}(a,d)\right)  \nonumber \\
&= n! \left(\sum_{i=1}^{n+1} \left(\frac{1}{a}\right)^i \frac{ \LLi_{s+i}(d)}{(n+1-i)!\,(m+n+2-i) } + \left(\frac{1}{a}\right)^{n+1}\LSI_{s+n+1}^{m}(a,d)\right).
\label{YequalsXCut}
\end{align}

Next, we consider the remaining cases where the interface $\Gamma$ is not parallel to one of the triangle edges.
For $a\ne0$, $b\ne0$ and $a+b\ne 0$
\begin{align}
\TRIB_s^{mn}&(a,b,d)=-\lim_{t \rightarrow \infty}\frac{1}{t^s}\int_{0}^{1} \int_0^x x^m y^n \Li_s(-\exp( (a x +b y +d )t ) \,dy\,dx  \nonumber \\
&= -\lim_{t \rightarrow \infty}\frac{1}{t^s} \int_0^1  n!
\left(\sum_{j=1}^{n+1} \frac{(-1)^{j-1}}{(n+1-j)!} x^{m+n+1-j} \frac{\Li_{s+j}(-\exp( ((a +b)x + d )t )}{(bt)^j} \right.\nonumber\\
&\qquad\qquad\qquad\qquad\qquad\qquad\quad\left.- 
(-1)^{n} x^m \frac{\Li_{s+n+1}(-\exp(a x + d )t )}{(bt)^{n+1}}
\right) dx \nonumber\\
&= -\sum_{j=1}^{n+1} \frac{n!}{(-b)^j(n+1-j)!} \LSI_{s+j}^{m+n+1-j}(a+b,d) + \frac{n!}{(-b)^{n+1}}\LSI_{s+n+1}^{m}(a,d).
\label{triangle}
\end{align}

Alternatively, the same integral could be evaluated by reversing the order of integration. Specifically,
\begin{align}
\TRIB_s^{mn}&(a,b,d)=-\lim_{t \rightarrow \infty}\frac{1}{t^s}\int_{0}^{1} \int_y^1 x^m y^n \Li_s(-\exp( (a x +b y +d )t ) \,dx\,dy  \nonumber \\
&= -\lim_{t \rightarrow \infty}\frac{1}{t^s} \int_0^1  m!
\left(\sum_{j=1}^{m+1} \frac{(-1)^{j-1}}{(m+1-j)!}y^n \frac{\Li_{s+j}(-\exp( by + a + d )t )}{(at)^j} \right.\nonumber\\
&\hspace{4cm}\left.- \frac{(-1)^{j-1}}{(m+1-j)!} y^{m+n+1-j} \frac{\Li_{s+j}(-\exp( ((a +b)y + d )t )}{(at)^j}
\right) dx \nonumber\\
&=\sum_{j=1}^{m+1} \frac{m!}{(m+1-j)!} \frac{-1}{(-a)}^j \left( \LSI_{s+j}^{m+n+1-j}(a+b,d) -\LSI_{s+j}^{n}(b,a+d)  \right).
\label{triangleB}
\end{align}
In the limit for $b\rightarrow0$, with $|a|>M>0$, Eq.~\eqref{triangle} may suffer from overflow. Similarly, in the limit for $a\rightarrow0$, with $|b|>M>0$, Eq.~\eqref{triangleB} may suffer from overflow. The choice of which formula to use, Eq.~\eqref{triangle} or Eq.~\eqref{triangleB}, should take into consideration the magnitude of $a$ and $b$.

\begin{remark} \label{rem:tri}
In Eq.~\eqref{triangle}, for $a+b+d\le0$ the summation within the $\LSI_{s+j}^{m+n+1-j}(a+b,d)$ terms vanishes. This is due to $\LLi_{s+1+i}(x) = 0$, with $s\ge-1$, $i\in \mathbb{Z}^+$, and non positive argument $x$. Specifically, for $a+b+d\le0$, Eqs.~\eqref{triangle} and \eqref{triangleB} reduce to
\begin{align}
{\TRIBR}_s^{mn}(a,b,d)=& n!\left( \frac{\LLi_{s+m+n+2}(d)}{(-(a+b))^{m+n+2}} 
\sum_{j=1}^{n+1} 
 \frac{(m+n+1-j)!}{(n+1-j)!}    \left(\frac{a+b}{b}\right)^j  \right. \nonumber \\
&\qquad +\left. \frac{\LSI_{s+n+1}^{m}(a,d)}{(-b)^{n+1}}\right),
\label{triangleReduced}
\end{align}
and
\begin{align}
{\TRIBR}_s^{mn}(a,b,d)&=
m!\left( - \frac{\LLi_{s+m+n+2}(d)}{(-(a+b))^{m+n+2}} 
\sum_{j=1}^{m+1} 
 \frac{(m+n+1-j)!}{(m+1-j)!}    \left(\frac{a+b}{a}\right)^j  \right. \nonumber \\
&\left. \qquad +  \frac{n!}{(-b)^{n+1}}
\sum_{j=1}^{m+1} 
 \frac{1}{(m+1-j)!} \frac{\LLi_{s+n+j+1}(a+d)}{(-a)^{j}}  \right),
\label{triangleBReduced}
\end{align}
which are less expensive to compute. 
For $s=-1$ and $a+b+d>0$, we can still take advantage of this reduction by changing the sign of the normal and utilizing Proposition \ref{stepFunction}.
Namely,
$$\TRIB_{-1}^{mn}(a,b,d) = \TRIBR_{-1}^{mn}(-a,-b,-d).$$
 A similar reasoning can be extended to the cases $\TRIB_s^{mn} (a,0,d)$, 
 $\TRIB_s^{mn} (0,b,d)$ and $\TRIB_s^{mn} (a,-a,0)$, when  $a+d\le0$,  $b+d\le0$, and  $d\le0$, respectively. However, special attention should be used if $s=-1$ and $a+b+d=0$, since for this case the first terms in the ``supposedly vanishing'' sums would be $\LLi_0(0) = -0.5 \ne 0$. Rewriting the three reduced formulas in a conservative way, always including the first term in the sum, leads to
 \begin{align}
&\TRIBR_s^{mn}(a,0,d)= \frac{1}{n+1}\left(-\frac{\LLi_{s+1}(a + d ) }{\,a} 
+ (m+n+1)! \frac{(-1)^{m+n+1} \LLi_{s+m+n+2}(d)}{a^{m+n+2}}\right),\label{trib0}\\
&\TRIBR_s^{mn}(0,b,d)= \frac{1}{m+1}\left( 
n!
\frac{(-1)^{n} \LLi_{s+n+1}(d)}{b^{n+1}} - 
 (m+n+1)! \frac{(-1)^{m+n+1} \LLi_{s+m+n+2}(d)}{b^{m+n+2}} \right),\label{tria0}\\
&\TRIBR_s^{mn}(a,-a,d)=  
\frac{\LLi_{s+1}(d )}{ (m+n+1) a}  
+ n! m! \sum_{i=1}^{m+1} \frac{(-1)^{i} \LLi_{s+n+1+i}(a+d ) }{(m+1-i)!\,a^{n+1+i}},  \label{tribma}
 \end{align}
 which hold for $a+b+d\le0$ and $s\ge-1$. 
\end{remark}
We also include the two alternative formulas below. These are obtained 
from \eugenio{the derivative identity in Eq.~\eqref{derivativeOfLi} and by recursive integration by parts},
increasing the monomial power and reducing the polylogarithm order $s$ until it reaches $-1$. Namely,
for $s\ge0$ and $b\ne0$, 
\begin{align} 
\TRIC_s^{m,n} (a,b,d)&=\dfrac{n!}{(n+s+1)!} (-b)^{s+1} \TRIA_{-1}^{m,n+s+1}(a,b,d) \nonumber \\
&+ \sum_{i=0}^{s} \dfrac{n!}{(n+i+1)!} (-b)^{i} \LSI_{s-i}^{m+n+i+1}(a+b,d), \label{triAlternative}
\end{align}
and, for $s\ge0$ and $a\ne 0$,

\begin{align} 
\TRIC_s^{m,n} (a,b,d)&=\dfrac{m!}{(m+s+1)!} (-a)^{s+1} \TRI_{-1}^{m+s+1,n}(a,b,d) \nonumber \\
&+ \sum_{i=0}^{s} \dfrac{m!}{(m+i+1)!} (-a)^{i} \left(\LSI_{s-i}^{n}(b,d+a)-\LSI_{s-i}^{m+n+i+1}(a+b,d)\right). \label{triAlternativeB}
\end{align}

For $a+b+d>\max(|a|,|b|)$, the combination of Remark \ref{rem:tri} and Eqs.~\eqref{triAlternative}-\eqref{triAlternativeB} yields a formulation which protects against overflow for $a \rightarrow 0$ and/or $b \rightarrow 0$. In particular, the calls to the $$\TRIA_{-1}^{m, n+s+1}(a,b,d)  \mbox{ and }  \TRIA_{-1}^{m+s+1,n}(a,b,d)$$ integrals in Eq.~\eqref{triAlternative} and Eq.~\eqref{triAlternativeB} can be replaced by 
 $$\TRIA_{-1}^{m, n+s+1}(-a,-b,-d)\mbox{ and }\TRIA_{-1}^{m+s+1,n}(-a,-b,-d),$$ respectively, for which $(-a) + (-b) +(-d) <0$.

At last we include the degenerate case when both $a=0$ and $b=0$ for $s \ge 0$, which was excluded because of the constraint $a^2+b^2>1$. This case is needed for external calls made by higher dimensional objects, such as the tetrahedron and prism, for which the normal $\bm n =\left<a,b,c \right>$ could take the form $\bm n = \left<0,0,c \right>$. After integration
\begin{equation} 
\TRIA_s^{mn}(0,0,d) \label{triDeg} =-\lim_{t \rightarrow \infty}\frac{1}{t^s}\int_{0}^{1} \int_0^x x^m y^n \Li_s(-\exp( d t ) ) dy\,dx = -\frac{\LLi_s^{mn}(d)}{(n+1)(m+n+2)}.
\end{equation}

The pseudo-code for the triangle integration is given in Algorithms \ref{algtri0} and \ref{algtri1}. Algorithm \ref{algtri0} evaluates the integral in Eq.~\eqref{TriOrigin} on the triangle
$\{(x,y):x\in[0,1],y\in[0,1-x]\}$. It calls the function {\sc Triangle\_Integration\_A}
in Algorithm \ref{algtri1}, which evaluates the transformed integral in Eq.~\eqref{TriAfterTransformation} on the triangle $\{(x,y):x\in[0,1],y\in[0,x]\}$.
{\sc Triangle\_Integration\_A} handles the degenerate case $a=b=0$ and sorts the different $s-$cases. For each case it ensures that the reduced integration function,  {\sc Triangle\_Integration\_BR}, is called only for $a+b+d\le0$. For $0<a+b+d\le\max(|a|,|b|)$, the 
function {\sc Triangle\_Integration\_B} is called, otherwise 
the alternative function {\sc Triangle\_Integration\_C} is used. The recursive calls follow from the patterns developed in Eqs.~\eqref{triAlternative} and \eqref{triAlternativeB}. Every time
the line segment integration formula,  $\LSI$, is needed the function {\sc Line\_Segment\_Integration} in Algorithm \ref{algLSI} is called.

\begin{algorithm} 
\caption{Pseudo-code for the integration of Eq.~\eqref{TriOrigin} on the triangle $\{(x,y):x\in[0,1],y\in[0,1-x]\}$ cut by the line $a\,x + b\,y+ d=0$ with $\bm n = \left<a,b \right>$, $\|\bm n\|>0$, $\bm m = \left<m,n \right>$, and $s=-1,0,1,\dots$. For $s=-1$ and $\|\bm n\|=1$ it corresponds to the interface integral. 
For $s=0$ it corresponds to the subdomain integral.}
\label{algtri0}
\begin{algorithmic}[1] 
  \Function {Triangle\_Integration}{$a$, $b$, $d$, $m$ $n$, $s$}
    \State{\Return $\text{\sc{Triangle\_Integration\_A}} (-a,\,b,\,d+a,\,m,\,n,\,s)$}
  \EndFunction
  \end{algorithmic}
\end{algorithm}

\begin{algorithm} 
\caption{Pseudo-code for the integration of Eq.~\eqref{TriAfterTransformation} on the triangle
$\{(x,y):x\in[0,1],y\in[0,x]\}$ cut by the line $a\,x + b\,y+ d=0$ with $\bm n = \left<a,b \right>$, $\bm m = \left<m,n \right>$ and $s=-1,0,1,\dots$. For $s=-1$ and $\|\bm n\|=1$ it corresponds to the interface integral. 
For $s=0$ and $\|\bm n\|>0$ it corresponds to the subdomain integral.}
\label{algtri1}
\begin{algorithmic}[1] 
  \Function {Triangle\_Integration\_A}{$a$, $b$, $d$, $m$ $n$, $s$}
  \vspace{0.5mm}
   \If{$b = 0$ and $a = 0$}
         \Return ${\TRI}_s^{mn}(0,0,d) $ from Eq.~\eqref{triDeg}
   \EndIf
  \vspace{0.5mm}
    \If{$s = -1$}
        \If{$a+b+d \le 0$} \Return  $\text{\sc{Triangle\_Integration\_BR}} (a,\,b,\,d,\,m,\,n,\,-1)$\vspace{0.5mm} 
        \Else{} \Return$\text{\sc{Triangle\_Integration\_BR}} (-a,\,-b,\,-d,\,m,\,n,\,-1)$\vspace{0.5mm} 
        \EndIf
    \Else
     \If{$a+b+d \le 0$} \Return  $\text{\sc{Triangle\_Integration\_BR}} (a,\,b,\,d,\,m,\,n,\,s)$ \vspace{0.5mm}
      \ElsIf{$a+b+d \le \max(|a|,|b|)$} \Return$\text{\sc{Triangle\_Integration\_B}} (a,\,b,\,d,\,m,\,n,\,s)$\vspace{0.5mm}
      \Else {} \Return$\text{\sc{Triangle\_Integration\_C}} (a,\,b,\,d,\,m,\,n,\,s)$
      \EndIf
    \EndIf
  \EndFunction
\end{algorithmic}
\vspace{0.25cm}

\begin{algorithmic}[1] 
  \Function {Triangle\_Integration\_B}{$a$, $b$, $d$, $m$ $n$, $s$}
    \If{$b = 0$}
         {\Return $\;\TRIB_s^{mn}(a,0,d) $ from Eq.~\eqref{verticalCut} \vspace{0.5mm}}
    \ElsIf{$a=0$}
      {\Return $\;\TRIB_s^{mn}(0,a,d) $ from Eq.~\eqref{horizontalCut}\vspace{0.5mm}} 
    \ElsIf{$a+b=0$} 
        {\Return $\;\TRIB^{mn}(a,-a,d) $ from Eq.~\eqref{YequalsXCut}\vspace{0.5mm}} 
    \Else  
        \If{$|a| \le |b|$}
            {\Return $\;\TRIB_s^{mn}(a,b,d) $ from Eq.~\eqref{triangle}\vspace{0.5mm}} 
        \Else {}  
            {\Return $\;\TRIB_s^{mn}(a,b,d) $ from Eq.~\eqref{triangleB}\vspace{0.5mm}}
        \EndIf
    \EndIf
  \EndFunction
\end{algorithmic}

\vspace{0.25cm}
\begin{algorithmic}[1] 
\Function {Triangle\_Integration\_BR}{$a$, $b$, $d$, $m$ $n$, $s$}
    \If{$b = 0$}
        {\Return $\;\TRIBR_s^{mn}(a,0,d) $ from Eq.~\eqref{trib0}\vspace{0.5mm}}
    \ElsIf{$a=0$}
      {\Return $\;\TRIBR_s^{mn}(0,a,d) $ from Eq.~\eqref{tria0}\vspace{0.5mm}} 
    \ElsIf{$a+b=0$} 
       {\Return $\;\TRIBR^{mn}(a,-a,d) $ from Eq.~\eqref{tribma}\vspace{0.5mm}} 
      
    \Else  
        \If{$|a| \le |b|$}
            {\Return $\;\TRIBR_s^{mn}(a,b,d) $ from Eq.~\eqref{triangleReduced} \vspace{0.5mm}} 
        \Else {}  
            {\Return $\;\TRIBR_s^{mn}(a,b,d) $ from Eq.~\eqref{triangleBReduced}\vspace{0.5mm}} 
        \EndIf
    \EndIf
  \EndFunction
  \end{algorithmic}
  \vspace{0.25cm}
  \begin{algorithmic}[1] 
\Function {Triangle\_Integration\_C}{$a$, $b$, $d$, $m$ $n$, $s$}
    \If{$|a|\le|b|$} \Return
        \State{ $ 
         \sum_{i=0}^{s} \dfrac{n!}{(n+i+1)!} (-b)^{i} \LSI_{s-i}^{m+n+i+1}(a+b,d)$} 
        \State{ $\qquad +\dfrac{n!}{(n+s+1)!} (-b)^{s+1}\text{\sc{Triangle\_Integration\_A}}(a,\,b,\,d,\,m,\,n+s+1,\,-1)  $}
    \Else {} \Return
        \State{ 
        $\sum_{i=0}^{s} \dfrac{m!}{(m+i+1)!} (-a)^{i} \left(\LSI_{s-i}^{n}(b,d+a)-\LSI_{s-i}^{m+n+i+1}(a+b,d)\right)$}
        \State{$\qquad  +\dfrac{m!}{(m+s+1)!} (-a)^{s+1} \text{\sc{Triangle\_Integration\_A}}(a,\,b,\,d,\,m+s+1,\,n,\,-1)$}
    \EndIf
    \EndFunction
  
  \end{algorithmic}
\end{algorithm}


\subsection{TTI: Tetrahedron Integration 
on $\left\{ \protect\begin{array}{@{}l@{}} 0\le x\le 1 \\ 0\le y\le 1-x \\0\le z \le 1-x-y
         \protect\end{array} \right.,$ with $a^2+b^2+c^2>0$}

To ease the computation in the case of the tetrahedron, we choose different polynomial bases depending on the magnitude of the coefficients $a$, $b$, and $c$. 

Let $m_1 = \max( |b-a| , |c-b|)$ and $m_2 = |a-c|$.
For $m_1 \ge m_2$, we evaluate integrals in the form
\begin{equation}
\TTI_s^{mno}(a,b,c,d) = -\lim_{t \rightarrow \infty}\frac{1}{t^s}\iiint_{\mbox{\sc{Tet}}}  (x+y+z)^m (y+z)^n z^o \Li_s(-\exp( (a x +b y + c z +d )t ) \,dV, \label{TET1}   
\end{equation}
else, we evaluate integrals in the form
\begin{equation}
\TTI_s^{mno}(a,b,c,d) = -\lim_{t \rightarrow \infty}\frac{1}{t^s} \iiint_{\mbox{\sc{Tet}}}    (y+z+x)^m (z+x)^n x^o   \Li_s(-\exp( (a x +b y + c z +d )t ) \,dV.  \label{TET2}     
\end{equation}

For $a = b = c$ (or $\max(m1,m2)=0$), integral \eqref{TET1} is considered and after integration we get
\begin{align}
\TTI_s^{mno}(a,a,a,d)= \frac{LSI_s^{m+n+o+2}(a,d)}{(o+1)(o+n+2)}. \label{TTIaeqbeqc}
\end{align}
Details of computation are given below. 

We make the following change of variables and constant renaming
\begin{itemize}
\item for Eq.~\eqref{TET1},
$$x'=x+y+z,\,y'=y+z,\,z'=z,\quad a'=a,\,b'=b-a,\, c'=c-b,\,d\,'=d,$$
\item for Eq.~\eqref{TET2}, 
$$x'=y+z+x,\,y'=z+x,\,z'=x,\quad a'=b,\,b'=c-b,\, c'=a-c,\,d\,'=d,$$
\end{itemize}
always obtaining the same integral
$$ -\lim_{t \rightarrow \infty}\frac{1}{t^s}\int_{0}^{1} \int_0^{x'} \int_0^{y'}   {x'}^{\,m} {y'}^{\,n} {z'}^{\,o} \Li_s(-\exp( (a' x' +b' y' + c' z' +d' )t ) \,dz'\,dy'\,dx'.$$
 
Dropping the $'$ superscript, for a fixed $s=-1,0,1,\dots$, the problem reduces to find integrals in the form
\begin{equation} \TTIA_s^{mno}(a,b,c,d)=-\lim_{t \rightarrow \infty}\frac{1}{t^s}\int_{0}^{1} \int_0^x \int_0^{y}   x^m y^n z^o \Li_s(-\exp( (a x +b y + c z +d )t ) \,dz\,dy\,dx\label{tetrahedronA}.
\end{equation}
The case $\max(|b|,|c|) = \max(m_1,m_2) = 0$ was already considered in Eq.\eqref{TTIaeqbeqc}.  
This corresponds to $\TTIA_s^{mno}(a,0,0,d)$, whose integration is straightforward.

Below, we consider only $\max(|b|,|c|) = \max(m_1,m_2)>0$.
For $|b| \le |c|$, after integrating in $z$
\begin{align}
\TTIB&_s^{mno}(a,b,c,d)\nonumber \\ &= -\lim_{t \rightarrow \infty}\frac{1}{t^s}\int_{0}^{1} \int_0^{x}    x^m y^n \Bigg(\sum_{i=1}^{o+1} \frac{o!(-1)^{i-1}}{(o+1-i)!} y^{o+1-i} \frac{\Li_{s+i}(-\exp( (a x +b y + cy + d )t )}{(ct)^i} \nonumber \\
&\quad - o! (-1)^{o}  \frac{\Li_{s+o+1}(-\exp( (a x +b y + d )t )}{(ct)^{o+1}} \Bigg)dy\,dx.\nonumber\\
&= -\lim_{t \rightarrow \infty}\frac{1}{t^s}\int_{0}^{1} \int_0^{x} \Bigg(\sum_{i=1}^{o+1} \frac{o!(-1)^{i-1}}{(o+1-i)!} x^{m} y^{n+o+1-i} \frac{\Li_{s+i}(-\exp( (a x + (b+c) y + d )t )}{(ct)^i} \nonumber \\
&\quad - o! (-1)^{o} x^m y^n  \frac{\Li_{s+o+1}(-\exp( (a x +b y + d )t )}{(ct)^{o+1}} \Bigg)dy\,dx.
\end{align}
Simplifying and using the triangle integration formula yields
\begin{align}
\TTIB_s^{mno}(a,b,c,d)=&  -\sum_{i=1}^{o+1} \frac{o!}{(o+1-i)!} \frac{1}{(-c)^i} \; 
\TRIA_{s+i}^{m,\;n+o+1-i}(a,b+c,d) + \frac{o!}{(-c)^{o+1}} \TRIA_{s+o+1}^{m n}(a,b,d). \label{TTIB1}
\end{align}
For $|c| < |b|$, we reverse the order of integration and after simplification get
\begin{align} \TTIB_s&^{mno}(a,b,c,d)=-\lim_{t \rightarrow \infty}\frac{1}{t^s}\int_{0}^{1} \int_0^x \int_z^{x}   x^m y^n z^o \Li_s(-\exp( (a x +b y + c z +d )t ) \,dy\,dz\,dx \nonumber
\\
& =  
\sum_{i=1}^{n+1} \frac{n!}{(n+1-i)!} \frac{1}{(-b)^i} \;\left( 
\TRIA_{s+i}^{m,\;n+o+1-i}(a,b+c,d) - \TRIA_{s+i}^{m+n+1-i,o}(a+b,c,d) \right). \label{TTIB2}
\end{align}

All limiting cases, are left to be handled by the triangle integration formula as described in the previous section.

For $s\ge0$ and $a +b +c+d > \max(|b|,|c|)$, we also include the alternative formulas below. These are obtained from \eugenio{the derivative identity in Eq.~\eqref{derivativeOfLi} and by recursive integration by parts},, increasing the monomial power and reducing the polylogarithm order $s$ until it reaches $-1$. Namely, For $|b|\le|c|$ we use
\begin{align} 
\TTIC_s^{m,n,o} (a,b,c,d) &=\dfrac{o!}{(o+s+1)!} (-c)^{s+1} \TTIA_{-1}^{m,n,o+s+1}(a,b,c,d) \nonumber \\
&+ \sum_{i=0}^{s} \dfrac{o!}{(o+i+1)!} (-c)^{i} \TRIA_{s-i}^{m, n+o+i+1}(a,b+c,d), 
\end{align}
otherwise 
\begin{align} 
\TTIC_s&^{m,n,o} (a,b,c,d) =\dfrac{n!}{(n+s+1)!} (-b)^{s+1} \TTIA_{-1}^{m,n+s+1,o}(a,b,c,d) \nonumber \\
&+ \sum_{i=0}^{s} \dfrac{n!}{(n+i+1)!} (-b)^{i}\left( \TRIA_{s-i}^{m+n+i+1,o}(a+b,c,d) - \TRIA_{s-i}^{m, n+o+i+1}(a,b+c,d) \right). 
\end{align}

The pseudo-code for the integration over the tetrahedron is given in Algorithms \ref{algtet0} and \ref{algtet1}.
Every time
the triangle integration formula  $\TRIA$ is needed, the function {\sc Triangle\_Integration\_A} in Algorithm \ref{algtri1} is called.

\begin{algorithm} 
\caption{Pseudo-code for the integration of Eqs.~\eqref{TET1}-\eqref{TET2} on the tetrahedron $\{(x,y,z):x\in[0,1],y\in[0,1-x], z\in[0,1-x-y]\}$ cut by  the plane $a\,x + b\,y+ c\,z + d=0$ with 
$\bm n = \left<a,b,c \right>$, $\|\bm n\|>0$,
$\bm m = \left<m,n,o \right>$ and $s=-1,0,1,\dots$. For $s=-1$ and $\|\bm n\|=1$ it corresponds to the interface integral. For $s=0$ it corresponds to the subdomain integral.}
\label{algtet0}
\begin{algorithmic}[1] 
  \Function {Tetrahedron\_Integration}{$a$, $b$, $c$, $d$, $m$ $n$, $o$, $s$}
  \If{$a=b=c$}
    \State{\Return $\TTI$ from Eq.~\eqref{TTIaeqbeqc}}
  \EndIf
  \State{$m_1 = \max( |a-b| , |c-b|)$, $m_2 = |a-c|$}
    \If{$m_1 \ge m_2$}
       \State{ \Return $\text{\sc{Tetrahedron\_Integration\_A}} (a,\,b-a,\,c-b,\,d,\,m,\,n,\,o,\,s)$}
    \Else
        \State{ \Return $\text{\sc{Tetrahedron\_Integration\_A}} (b,\,c-b,\,a-c,\,d,\,m,\,n,\,o,\,s)$}
    \EndIf
  \EndFunction
  \end{algorithmic}
\end{algorithm}
\begin{algorithm} 
\caption{Pseudo-code for the integration of Eq.~\eqref{tetrahedronA} on the tetrahedron $\{(x,y,z):x\in[0,1],y\in[0,x], z\in[0,y]\}$ cut by the plane $a\,x + b\,y+ c\,z + d=0$ with 
$\bm n = \left<a,b,c \right>$, either $b\ne0$ or $c \ne 0$, $\bm m = \left<m,n,o \right>$ and $s=-1,0,1,\dots$. For $s=-1$ and $\|\bm n\|=1$ it corresponds to the interface integral. For $s=0$ it corresponds to the subdomain integral.
}

\label{algtet1}
\begin{algorithmic}[1] 
  \Function {Tetrahedron\_Integration\_A}{$a$, $b$, $c$, $d$, $m$ $n$, $o$, $s$}
    \State{$sum = a + b + c + d$}
    \If{$s = -1$}
        \If {$sum \le 0$} 
            \State {\Return $\text{\sc{Tetrahedron\_Integration\_B}} (a,\,b,\,c,\,d,\,m,\,n,\,o,\,-1)$}
        \Else
            \State{\Return $\text{\sc{Tetrahedron\_Integration\_B}} (-a,\,-b,\,-c,\,-d,\,m,\,n,\,o,\,-1)$}
        \EndIf
    \Else
        \If{ $sum \le \max(|b|,|c|)$}
            \State{\Return $\text{\sc{Tetrahedron\_Integration\_B}} (a,\,b,\,c,\,d,\,m,\,n,\,o,\, s)$}
        \Else
            \State{\Return $\text{\sc{Tetrahedron\_Integration\_C}} (a,\,b,\,c,\,d,\,m,\,n,\,o,\, s)$}
        \EndIf
  \EndIf
  
  \EndFunction
  \end{algorithmic}
  \vspace{0.25cm}

\begin{algorithmic}[1] 
  \Function {Tetrahedron\_Integration\_B}{$a$, $b$, $c$, $d$, $m$ $n$, $o$, $s$}
        \vspace{1mm}
        \If{$|b|\le|c|$}  
        \Return $\TTIB$ from Eq.~\eqref{TTIB1} \vspace{1mm}
        \Else {} \Return $\TTIB$ from Eq.~\eqref{TTIB2} \vspace{1mm}
        \EndIf
  \EndFunction
  \end{algorithmic}
  \vspace{0.25cm}
  \begin{algorithmic}[1] 
  \Function {Tetrahedron\_Integration\_C}{$dim$, $\bm n$, $d$, $\bm m$, $s$}
    \If{$|b|\le|c|$} \Return
        \State{ $ 
         \sum_{i=0}^{s} \dfrac{o!}{(o+i+1)!} (-c)^{i} \TRIA_{s-i}^{m, n+o+i+1}(a,b+c,d)$} 
        \State{ $\qquad +\dfrac{o!}{(o+s+1)!} (-c)^{s+1}\text{\sc{Tetrahedron\_Integration\_A}}(a,\,b,\,c,\,d,\,m,\,n,\,o+s+1,\,-1)  $}
    \Else {} \Return
        \State{ 
        $ \sum_{i=0}^{s} \dfrac{n!}{(n+i+1)!} (-b)^{i}\left( \TRIA_{s-i}^{m+n+i+1,o}(a+b,c,d) - \TRIA_{s-i}^{m, n+o+i+1}(a,b+c,d) \right)$}
        \State{$\qquad  +\dfrac{n!}{(n+s+1)!} (-b)^{s+1} \text{\sc{Tetrahedron\_Integration\_A}}(a,\,b,\,c,\,d,\,m,\,n+s+1,\,o,\,-1)$}
    \EndIf
  \EndFunction
  \end{algorithmic}
\end{algorithm} 

\subsection{PRI: Prism Integration 
on 
$\left\{ \protect\begin{array}{@{}l@{}} 0\le x\le 1 \\ 0\le y\le 1-x \\-1\le z \le 1
         \protect\end{array} \right.,$ with $a^2+b^2+c^2>0$}

The implementation of a polynomial basis, whose elements are given by $(1-x)^m y^n z^o$, allows for  computational simplicity when considering integrals in the form
\begin{equation}
\PRI = -\lim_{t \rightarrow \infty}\frac{1}{t^s}\iiint_{\mbox{\sc{Pri}}} (1-x)^m y^n \left(\frac{1+z}{2}\right)^o \Li_s(-\exp( (a x +b y+c z +d )t ) \frac{dV}{2}.    \label{PriOrigin}
\end{equation}
By using the following transformation 
$$ x' = 1-x, \, y' = y, \, z' = \frac{1+z}{2}, \quad a'=-a, b'= b, c'=2 c,\, d' = d + a - c$$
we obtain
$$ -\lim_{t \rightarrow \infty}\frac{1}{t^s}\int_{0}^{1} \int_0^{x'} \int_{0}^1 {x'}^{\,m} {y'}^{\,n} {z'}^{\,o} \Li_s(-\exp( (a' x' +b' y' + c' z' + d' )t )) dz' \,dy'\,dx'$$

Dropping the $'$ superscript, for a fixed $s=-1,0,1,\dots$, the problem reduces to integrals in the form
$$ \PRIA_s^{mno}(a,b,c,d)=-\lim_{t \rightarrow \infty}\frac{1}{t^s}\int_{0}^{1} \int_0^x \int_{0}^1   x^m y^n z^o \Li_s(-\exp( (a x +b y + c z +d )t ) \,dz\,dy\,dx.$$ 

For $|c|\ge\max(|a|,|b|)$, after integrating in z we get
\begin{align}
\PRIB&_s^{mno}(a,b,c,d)
= -\lim_{t\rightarrow \infty}  \int_{0}^{1}\int_{0}^{x} x^m y^n \left(- 
\sum_{i=1}^{o+1} \frac{o!} {(o+1-i)!} \, \frac{1}{(-c)^i} \frac{\Li_{s+i} (-\exp( (a x+by+c +d )t )}{ t^{s+i}} \nonumber\right. \\
&  
\left. \qquad + \frac{o!}{(-c)^{o+1}} \frac{\Li_{s+o+1} (-\exp( (ax +by + d )t )}{ t^{s+k}} \right) \, dy \, dx \nonumber\\
=& -\sum_{i=1}^{o+1} \frac{o!} {(o+1-i)!} \, \frac{1}{(-c)^i}   \TRIA_{s+i}^{mn}(a,b,c+d) +  \frac{o!}{(-c)^{o+1}} \TRIA_{s+o+1}^{mn}(a,b,d) . \label{PRIB1}
\end{align}
For $|b|\ge|a|$, after integrating first in $y$ and simplifying we have
\begin{align}
\PRIB&_s^{mno}(a,b,c,d)
= -\lim_{t \rightarrow \infty}\frac{1}{t^s}\int_{0}^{1} \int_0^1 \int_{0}^x   x^m y^n z^o \Li_s(-\exp( (a x +b y + c z +d )t ) \,dy\,dz\,dx \nonumber\\
=& -\sum_{i=1}^{n+1} \frac{n!} {(n+1-i)!} \, \frac{1}{(-b)^i}   \HCI_{s+i,2}^{\langle m+n+1-i,o \rangle}(\langle a+b,c \rangle,d) +  \frac{n!}{(-b)^{n+1}} \HCI_{s+n+1,2}^{\langle m,o \rangle}(\langle a,c\rangle, d) . \label{PRIB2}
\end{align}
Lastly, for all other cases,  after integrating first in $x$ and simplifying we obtain
\begin{align}
\PRIB&_s^{mno}(a,b,c,d)
= -\lim_{t \rightarrow \infty}\frac{1}{t^s}\int_{0}^{1} \int_0^1 \int_{y}^1   x^m y^n z^o \Li_s(-\exp( (a x +b y + c z +d )t ) \,dx\,dy\,dz \nonumber\\
=& \sum_{i=1}^{m+1} \frac{m!} {(m+1-i)!} \, \frac{1}{(-a)^i}   \left(-\HCI_{s+i,2}^{\langle n,o \rangle}(\langle b,c \rangle,a+d) +  \HCI_{s+i,2}^{\langle m+n+1-i,o \rangle}(\langle a+b,c\rangle, d)\right). \label{PRIB3}
\end{align}
All limiting cases are left to be handled by the triangle and the hypercube integration formulas previously described.

For $s\ge0$ and $a +b +c+d > \max(|a|,|b|,|c|)$, we also include the alternative formulas below. These are obtained from \eugenio{the derivative identity in Eq.~\eqref{derivativeOfLi} and by recursive integration by parts}, increasing the monomial power and reducing the polylogarithm order $s$ until it reaches $-1$. Namely, for $|c|\ge\max(|a|,|b|)$, we utilize the formula
\begin{align} 
\PRIC_s&^{m,n,o} (a,b,c,d) =\dfrac{o!}{(o+s+1)!} (-c)^{s+1} \PRIA_{-1}^{m,n,o+s+1}(a,b,c,d) \nonumber \\
&+ \sum_{i=0}^{s} \dfrac{o!}{(o+i+1)!} (-c)^{i} \TRIA_{s-i}^{m, n}(a,b,c+d), 
\end{align} and for $|b|\ge|a|$ we implement
\begin{align} 
\PRIC_s&^{m,n,o} (a,b,c,d) =\dfrac{n!}{(n+s+1)!} (-b)^{s+1} \PRIA_{-1}^{m,n+s+1,o}(a,b,c,d) \nonumber \\
&+ \sum_{i=0}^{s} \dfrac{n!}{(n+i+1)!} (-b)^{i} 
\HCI_{s-i,2}^{\langle m+n+i+1,o \rangle}(\langle a+b,c\rangle, d). 
\end{align}
For any other case we employ 
\begin{align} 
\PRIC_s&^{m,n,o} (a,b,c,d) =\dfrac{m!}{(m+s+1)!} (-b)^{s+1} \PRIA_{-1}^{m+s+1,n,o}(a,b,c,d) \nonumber \\
&+ \sum_{i=0}^{s} \dfrac{m!}{(m+i+1)!} (-a)^{i}\left( \HCI_{s-i,2}^{\langle n,o \rangle}(\langle b,c\rangle, a + d) -
\HCI_{s-i,2}^{\langle m+n+i+1,o \rangle}(\langle a+b,c\rangle, d) \right).
\end{align}

The pseudo-code for integration over the prism is given in Algorithms \ref{algpri0} and \ref{algpri1}.
Every time
the triangle integration formula  $\TRIA$ and the hypercube integration formula $\HCI$ are used, the functions {\sc Triangle\_Integration\_A} in Algorithm \ref{algtri1} and 
{\sc HyperCube\_Integration} in Algorithm \ref{alg1} are called.
 
 \begin{algorithm} 
\caption{Pseudo-code for the integration of Eq.~\eqref{PriOrigin} on the prism $\{(x,y,z):x\in[0,1],y\in[0,1-x], z\in[-1,1]\}$ cut by the plane $a\,x + b\,y+ c\,z + d=0$ with 
$\bm n = \left<a,b,c \right>$, $\|\bm n\|>0$,
$\bm m = \left<m,n,o \right>$ and $s=-1,0,1,\dots$. For $s=-1$ and $\|\bm n\|=1$ it corresponds to the interface integral. For $s=0$ it corresponds to the subdomain integral. 
}
\label{algpri0}
\begin{algorithmic}[1] 
  \Function {Prism\_Integration}{$a$, $b$, $c$, $d$, $m$ $n$, $o$, $s$}
    \State{\Return $\text{\sc{Prism\_Integration\_A}} (-a,\,b,\,2c,\,d+a-c,\,m,\,n,\,o,\,s)$}
  \EndFunction
  \end{algorithmic}
\end{algorithm}

\begin{algorithm} 
\caption{Pseudo-code for the integration on the prism $\{(x,y,z):x\in[0,1],y\in[0,x], z\in[0,1]\}$ cut by the plane $a\,x + b\,y+ c\,z + d=0$ with 
$\bm n = \left<a,b,c \right>$, $\|\bm n\|>0$,
$\bm m = \left<m,n,o \right>$ and $s=-1,0,1,\dots$. For $s=-1$ and $\|\bm n\|=1$ it corresponds to the interface integral. For $s=0$ it corresponds to the subdomain integral. }
\label{algpri1}
\begin{algorithmic}[1] 
  \Function {Prism\_Integration\_A}{$a$, $b$, $c$, $d$, $m$ $n$, $o$, $s$}
    \State{$sum = a + b + c + d$}
    \If{$s = -1$}
        \If {$sum \le 0$} 
            \State {\Return $\text{\sc{Prism\_Integration\_B}} (a,\,b,\,c,\,d,\,m,\,n,\,o,\,-1)$}
        \Else
            \State{\Return $\text{\sc{Prism\_Integration\_B}} (-a,\,-b,\,-c,\,-d,\,m,\,n,\,o,\,-1)$}
        \EndIf
    \Else
        \If{ $sum \le \max(|a|,|b|,|c|)$}
            \State{\Return $\text{\sc{Prism\_Integration\_B}} (a,\,b,\,c,\,d,\,m,\,n,\,o,\, s)$}
        \Else
            \State{\Return $\text{\sc{Prism\_Integration\_C}} (a,\,b,\,c,\,d,\,m,\,n,\,o,\, s)$}
        \EndIf
  \EndIf
  \EndFunction
  \end{algorithmic}
  
   \vspace{0.25cm}
\begin{algorithmic}[1] 
  \Function {Prism\_Integration\_B}{$a$, $b$, $c$, $d$, $m$ $n$, $o$, $s$}
        \vspace{1mm}
        \If{$|c|\ge\max(|a|,|b|)$}  
        \Return $\PRIB$ from Eq.~\eqref{PRIB1} \vspace{1mm}
        \ElsIf {$|b|>|a|$} \Return $\PRIB$ from Eq.~\eqref{PRIB2} \vspace{1mm}
        \Else{} \Return $\PRIB$ from Eq.~\eqref{PRIB3} \vspace{1mm}
        \EndIf
  \EndFunction
  \end{algorithmic}
  \vspace{0.25cm}
  \begin{algorithmic}[1] 
  \Function {Prism\_Integration\_C}{$a$, $b$, $c$, $d$, $m$ $n$, $o$, $s$}
    \If{$|c|\ge\max(|a|,|b|)$} \Return
        \State{ $ 
         \sum_{i=0}^{s} \dfrac{o!}{(o+i+1)!} (-c)^{i} \TRIA_{s-i}^{m, n}(a,b,c+d)$} 
        \State{ $\qquad +
        \dfrac{o!}{(o+s+1)!} (-c)^{s+1}
        \text{\sc{Prism\_Integration\_A}}(a,\,b,\,c,\,d,\,m,\,n,\,o+s+1,\,-1)$}
    \ElsIf {$|b|\ge|a|$} \Return
        \State{ 
        $\sum_{i=0}^{s} \dfrac{n!}{(n+i+1)!} (-b)^{i} 
\HCI_{s-i,2}^{\langle m+n+i+1,o \rangle}(\langle a+b,c\rangle, d) $}
        \State{$\qquad  +\dfrac{n!}{(n+s+1)!} (-b)^{s+1} \text{\sc{Prism\_Integration\_A}}(a,\,b,\,c,\,d,\,m,\,n+s+1,\,o,\,-1)$}
    \Else{} \Return    
    \State{ 
        $\sum_{i=0}^{s} \dfrac{m!}{(m+i+1)!} (-a)^{i}\left( \HCI_{s-i,2}^{\langle n,o \rangle}(\langle b,c\rangle, a + d) -
\HCI_{s-i,2}^{\langle m+n+i+1,o \rangle}(\langle a+b,c\rangle, d) \right)$}
        \State{$\qquad  +\dfrac{m!}{(m+s+1)!} (-a)^{s+1} \text{\sc{Prism\_Integration\_A}}(a,\,b,\,c,\,d,\,m+s+1,\,n,\,o,\,-1)$}
    \EndIf
  \EndFunction
  \end{algorithmic}
  
\end{algorithm}

\section{Note on the equivalent polynomial}
The equivalent polynomial problem can be stated as follows: find the equivalent polynomial coefficients $\bm c_o$, such that $M \bm c_o = \bm f_o$,  where
  \begin{align*}
    \bm f_o &= -\lim_{t\rightarrow \infty} t^{-s} \begin{pmatrix}
           \int_{\Omega} \bm b_{o,0}(\bm x) \Li_s (-\exp((\bm n \cdot \bm x+d)t)) \,d \bm x \\ 
            \\
           \int_{\Omega} \bm b_{o,1}(\bm x) \Li_s (-\exp((\bm n \cdot \bm x+d)t)) \,d \bm x \\
           \vdots \\
          \int_{\Omega} \bm b_{o,L}(\bm x) \Li_s (-\exp((\bm n \cdot \bm x+d)t)) \,d \bm x
         \end{pmatrix}
         \end{align*}
         and
          \begin{align*}
          M = 
\begin{pmatrix}
\int_{\Omega} \bm b_{o,0}(\bm x) \, \bm b_{o,0}(\bm x)\,d \bm x & \int_{\Omega} \bm b_{o,1}(\bm x) \, \bm b_{o,0}(\bm x) \,d \bm x & \cdots & \int_{\Omega} \bm b_{o,L}(\bm x)\, \bm b_{o,0}(\bm x)\,d \bm x \\ 
  \\
\int_{\Omega} \bm b_{o,0}(\bm x) \, \bm b_{o,1}(\bm x)\,d \bm x & \int_{\Omega} \bm b_{o,1}(\bm x) \, \bm b_{o,1}(\bm x) \,d \bm x & \cdots & \int_{\Omega} \bm b_{o,L}(\bm x) \, \bm b_{o,1}(\bm x)\,d \bm x \\
\vdots  & \vdots  & \ddots & \vdots  \\
\int_{\Omega} \bm b_{o,0}(\bm x) \, \bm b_{o,L}(\bm x)\,d \bm x & \int_{\Omega} \bm b_{o,1}(\bm x) \, \bm b_{o,L}(\bm x) \,d \bm x & \cdots & \int_{\Omega} \bm b_{o,L}(\bm x) \, \bm b_{o,L}(\bm x)\,d \bm x 
\end{pmatrix},
\end{align*}
with $s=-1$ or $0$. Here ${\bm b}_o$ is the basis of the polynomial space. Then, the equivalent polynomial is given by $p(\bm x) = \bm c_{o}^T \cdot {\bm b}_o.$
In order to avoid an ill-conditioned Gram matrix $M$, we implement orthogonal polynomials, via Grahm-Schmidt orthogonalization, using the $L^2$ inner product\cite{saad1986condition,xu1993multivariate}. This yields the following relation for basis elements: $\bm b_n = A \bm b_o$, where the components in the new basis, $\bm b_n$, are a linear combination of the components in the old basis, $\bm b_o.$ 
The matrix $A$ is an $(L+1) \times (L+1)$ lower triangular matrix, where $L$ is the dimension of the space spanned by the basis vector $\bm b_o$. 

The implementation of equivalent polynomial using an orthonormal basis yields  $$I \bm{c}_n = \bm{f}_n = A \bm{f}_o,$$ 
resulting in $$ p(\bm x) = (\bm{c}_n)^T \bm{b}_n  = \bm{f}_o^T A^T A \bm{b}_o(\bm x). $$
Note that the term $A^T A \bm{b}_o(\bm x)$ is independent of the hyperplane cut and can be evaluated off-line. Instead $\bm f_o$ changes and has to be recalculated for every new cut.

To this end, the continuous dependence of $\bm{f}_o$ with respect to the coefficients of the cut planes $\bm{n}$ and $d$ is of great help. Namely, for each considered element, we can explicitly evaluate 
$$ \frac{\partial \bm{f}_o}{\partial \bm{n}} \quad \mbox{ and } \quad  \frac{\partial \bm{f}_o}{\partial d},$$ and prove differentiability almost everywhere of $\bm f_0$ with respect to these parameters.
This implies that for each element a given set of quadrature rules can be evaluated and stored off-line, and a new quadrature integration rule can be reconstructed on-line by interpolation at very little cost and to any accuracy, making this technique far superior to any other existing method.

In 2D for a given line $a x + b y + d = 0$ we use the two parameter family given by the polar angle $\theta = \mbox{atan2}(b, a)$ and the $x-$intercept between the given line  and the lines
\begin{align*} 
& x = y , &&\mbox{ if $\theta$ is in the $1^{st}$ or $3^{rd}$ quadrant, or}\\
&1 - x= y, &&\mbox{ if $\theta$ is in the $2^{nd}$ or $4^{th}$ quadrant,}
\end{align*} respectively.

In 3D for a given plane $a x + b y +c z + d = 0$ we use the three parameter family given by the polar angle $\theta = \mbox{atan2}(b, a)$, the azimuthal angle $\phi = \mbox{acos}(c / \sqrt{a^2+b^2+c^2})$, 
and the $x-$intercept between the given plane and the lines
\begin{align*}
& x = y = z, && \mbox{ if $\theta$ and $\phi$ are in the $1^{st}$ or $7^{th}$ octant, or}\\
&  1-x = y = z, && \mbox{ if $\theta$ and $\phi$ are in the $2^{nd}$ or $8^{th}$ octant, or} \\
&  x = y = 1-z, && \mbox{ if $\theta$ and $\phi$ are in the $3^{rd}$ or $5^{th}$ octant, or} \\
&  x = 1 - y = z, && \mbox{ if $\theta$ and $\phi$ are in the $4^{th}$ or $6^{th}$ octant,} 
\end{align*}
respectively. 

These choices assure that for lines or planes cutting any of the considered elements the $x-$intercept is always in the interval $[0,1]$. Then, for each quadrant or octant, as above, we construct off-line matrices of coefficients spanning the whole range of parameters, and use on-line Lagrange interpolation to reconstruct the values of the coefficients for any $(x, \theta) \in [0,1] \times [-\pi, \pi]$ in 2D, or 
$(x, \theta, \phi) \in [0,1] \times [-\pi, \pi] \times [0,\pi]$ in 3D. 

Rather than storing and interpolating the equivalent polynomial coefficients, we store and interpolate the values of the equivalent polynomial evaluated at the quadrature points, for a given quadrature rule. This becomes particularly useful in the case of the tetrahedron where we used two different bases in the parent element depending of the values of the normal $\mathbf{n}$. While the coefficients $\mathbf{c}_n$ would differ for the two bases, the equivalent polynomial $p(x) = \mathbf{c}_n \mathbf{b}_n$ remains the same. Thus interpolation is still possible even when using interpolants evaluated with different bases.

\section{Conclusion}
The many closed-form algebraic expressions provided in the current work can easily be implemented into numerous PDE solvers when discontinuous functions are implemented. We have eliminated the need to consider complicated subdomains while simultaneously eliminating any error produced by a regularization parameter and polylogarithm approximation. We provide exact formulas for cumbersome subdomain and interface integrals, along with the associated algorithms. These closed-forms were designed with floating point arithmetic in mind. The results of this work provide one with the tools to eliminate many of the problems posed by discontinuous function integration. In this work, the discontinuities we considered were points, lines, and planes. \eugenio{Analytical integration on subdomains bounded by curved surfaces is currently being investigated. We have shown that analytical integration is still possible for elements cut by surfaces as complex as $$ P_n(x) + y ( a x + b ) + c z + d = 0,$$ for any degree $P_n$ polynomial. A preliminary version of this result is already available in the PhD thesis of the second author \cite{loftin2022exact} and will be analyzed in details in a forthcoming paper.}

\section*{Funding}
This work was supported by the National Science Foundation (NSF) Division of Mathematical Sciences (DMS) program, project 1912902.\\
The authors have no conflicts of interest to declare that are relevant to the content of this article.

\section*{Data availability}
Data sharing is not applicable to this article as no datasets were generated during the current study.

\eugenio{
\section*{Appendix A} 
\begin{lemma}
Let $D$ be a bounded connected domain with smooth boundary $\partial D$. Let $G(\bm x)$ be a smooth level set function. Let $\Gamma = \left\{ \bm x\in D : G(\bm x) = 0 \right\}$ be a continuous smooth embedded interface, that separates $D$ in the two subregions $D_1$ and $D_2$, such that $G(\bm x)>0 $ for all $ \bm x \in D_1$ and $G(\bm x)<0$ for all $ \bm x \in D_2$. Assume the measure $\mu(\Gamma \cap \partial D) = 0$. Then, for any differentiable function $f(\bm x)$
\begin{align*}
-\lim_{t\rightarrow\infty}t& \int_D f  \Li_{-1}(-\exp(G\,t)) \| \nabla G \| \,d\bm x
= \int_D f  \, \delta(G) \|\nabla(G)\| \,d\bm x.
\end{align*}
\end{lemma}
}
\eugenio{The $\|\nabla G\|$ term in both sides is needed since the level set $G(\bm x)$ only approximates the required condition, $\|\nabla d\|=1$, for a true distance $d(\bm x)$, see Appendix in \cite{holdych2008quadrature}. }

\begin{proof}
\eugenio{
Let $\bm n$ on $\partial D$ be the outer unit normal vector to $D$. 
Let $\widehat{\bm n} =-\frac{\nabla G}{\| \nabla G \|}$ be defined everywhere on $D$.
$\widehat{\bm n}$ is the unit vector orthogonal to the the level curves $G(\bm x) = const$, pointing in the direction of maximum decrease.
On the interface $\Gamma$, $\widehat{\bm n}$ is the unit outer normal to $D_1$.
Let $\partial D_1 = \partial D \cap D_1.$ 
Then, the boundary of $D_1$ is piece-wise-defined by $\partial D_1 \cup \Gamma$, with outer unit normal vectors $\bm n$ and $\widehat{\bm n}$, respectively.}

\eugenio{
Observe that by using the chain rule and the derivative property of the polylogarithm, we have
\begin{align} 
\nabla&\Li_0(-\exp(G\,t)) \cdot \widehat{\bm n} = \nabla\Li_0(-\exp(G\,t)) \cdot \left(-\frac{\nabla G}{\| \nabla G \|} \right) \nonumber \\
&= -t \frac{d}{d G} \Li_0(-\exp(G\,t)) \nabla G \cdot\frac{\nabla G}{\| \nabla G \|} = -t \Li_{-1}(-\exp(G\,t)) \| \nabla G \|. \label{derivative2}
\end{align}
Then,
\begin{align*}
-&\lim_{t\rightarrow\infty}t \int_D f\, \Li_{-1}(-\exp(G\,t)) \| \nabla G \| \,d\bm x=
\lim_{t\rightarrow\infty}\int_D \nabla\Li_0(-\exp(G\,t)) \cdot \left( f\, \widehat{\bm n}\right) \,d\bm x & \text{by Eq.~\eqref{derivative2}}
\\&= -\lim_{t\rightarrow\infty}\left(\int_D  \Li_0(-\exp(G\,t)) \nabla \cdot 
\left(f \,\widehat{\bm n} \right) \,d\bm x + \int_{\partial D}   \Li_0(-\exp(G\,t)) \, f \, \widehat{\bm n} \cdot \bm n \, dS \right) & \text{by Div. Thm.} 
\\& = \int_D  \U(G) \,\nabla \cdot \left(f \,\widehat{\bm n}\right) \,d\bm x - \int_{\partial D}   \U(G)\,f \, \widehat{\bm n}\cdot \bm n \, dS & \text{ by Eq.~\eqref{heaviside_equality}}
\\& = \int_{D_1}  \nabla \cdot \left(f \,\widehat{\bm n}\right) \,d\bm x - \int_{\partial D_1}  f \, \widehat{\bm n} \cdot \bm n \, dS & \text{ by Eq.~\eqref{heaviside_hmax}}\\
&= \left(\int_{\partial D_1}  f \,\widehat{\bm n}\cdot \bm n \,d S + 
\int_{\Gamma}  f \,\widehat{\bm n}\cdot \widehat{\bm n} \,d S \right) - \int_{\partial D_1}  f \, \widehat{\bm n} \cdot \bm n \, dS & \text{by Div. Thm.}
\\&= \int_{\Gamma}  f dS = \int_D f  \, \delta(G) \|\nabla G\| \,d\bm x.& \text{by Dirac delta Def.}
\end{align*}
Note that the proof holds only if the measure $\mu(\Gamma \cap \partial D) = 0$, for an appropriate product measure $\mu $, since from the third line to the fourth line, the integral equality on the boundary
$$\int_{\partial D}   \U(G)\,f \, \widehat{\bm n}\cdot \bm n \, dS =
\int_{\partial D_1}  f \, \widehat{\bm n} \cdot \bm n \, dS $$
is true only if the Heaviside function $\U$ is almost everywhere $1$ on $\partial D_1$ and almost everywhere $0$ on its complement. For $\mu(\Gamma \cap \partial D) \ne 0$, we would have measurable parts of the boundary $\partial D$ with $\U =0.5$, and the equality would not hold.}
\end{proof}

\bibliographystyle{plain}


\end{document}